# Concurrent multi-domain simulations in linear structural dynamics using multiple grid multiple time-scale (MGMT) method.

1)  Tejas H. Ruparel (corresponding author)
    School of Engineering and Applied Science,
    The George Washington University,
    Washington, DC 20052, USA
    Phone: 703-726-8323
    E-mail: truparel@gwu.edu

2)  Azim Eskandarian
    School of Engineering and Applied Science,
    The George Washington University,
    Washington, DC 20052, USA
    E-mail: eska@gwu.edu

3)  James D. Lee
    School of Engineering and Applied Science,
    The George Washington University,
    Washington, DC 20052, USA
    E-mail: truparel@gwu.edu

**Abstract**
Work presented in this paper describes a general algorithm and its finite element implementation for performing concurrent multiple sub-domain simulations in linear structural dynamics. Using this approach one can solve problems in which the domain under analysis can be selectively discretized spatially and temporally, hence allowing the user to obtain a desired level of accuracy in critical regions while improving computational efficiency globally. The mathematical background for this approach is largely based upon the fundamental principles of domain decomposition methods (DDM) and Lagrange multipliers, used to obtain coupled equations of motion for distinct regions of a continuous domain. These methods when combined together systematically yield constraint forces that not only ensure conservation of energy but also enforce continuity of velocities across sub-domain interfaces. Multiple grid (MG) connections between non-conforming sub-domains are modeled using mortar elements whereas coupled multiple time-scale (MT) equations are derived for the classical Newmark integration algorithm (and its constituents). Fully discretized equations of motion for component sub-domains, augmented with an interface condition are then solved using block elimination method and Crout factorization.

A proof of stability is provided using Energy method and overall efficiency, accuracy and stability of multiple sub-domain coupling is evaluated using a series of numerical examples. Primary observations are made for the evolution and distribution of kinematic quantities and structural wave propagation across connecting sub-domain. Numerical stability is verified by ensuring energy balance at global as well as component sub-domain level. Discussed examples highlight the greatest advantage of MGMT method; which is high simulation speedups (at the cost of reasonably small errors).

## 1. Introduction

Space and time are inherently coupled in the analysis of engineering problems. An approximate solution to these problems is usually obtained by using a mathematical tool suitable for the scale at which the physical phenomenon is addressed. Finite element method (FEM) is one such numerical technique that is used for solving problems in structural dynamics. It discretizes the governing equations over finite space and time and uses variational principles to minimize an error function in order to produce a stable solution. Numerical techniques as such, innately introduce discretization error due to the choice of finite space and time resolution [1,2]. Overall quality of the solution can be improved by spatial and temporal refinements, however at the expense of increasing number of unknowns and consequently longer computation times.

Over the past couple of decades, researchers have devoted a significant amount of effort towards modeling and implementing structural behavior at different resolutions in order to improve overall quality of the solution whilst preserving computational efficiency. Static mesh transition, adaptive mesh refinement (AMR) [3], mortar finite element methods [4,5] and finite element tearing and interconnecting (FETI) [6] are among the few techniques that are widely used to segregate problem size (total number of unknowns) between critical and remote regions. These techniques not only provide complete control over grid resolution, compared to fixed coarse or fine scale discretization, but also help in capturing local gradients and wave dynamics more accurately, hence yielding a better solution within the qualified region of interest. Domain decomposition method (DDM) [7,8] is another approach that enables effective implementation of these techniques by decomposing the domain under consideration into component sub-domains, which can then be modeled independently. Within the context of time discretization, a usual approach in coupling sub-domains is to use the same integration method (implicit or explicit) with the same time-step ($\Delta t$) globally over multiple grids. This is not recommended since it restricts one to analyze an entire domain using a single time-step that meets the stability and accuracy criteria for all elements. This is not desirable in the case of large scale problems since different regions could very well represent different stability and accuracy requirements. In addition, different regions may exhibit high frequency (wave propagation type) and/or low frequency (vibration type) responses, requiring explicit and/or implicit methods respectively. Accordingly it is much more economical to use different time-steps or different time-stepping algorithms in different sub-domains in order to capture local behavior as accurately as possible.

Mixed methods [9–13] are one classical approach where the time marching algorithm uses the same time-step but different schemes (implicit/explicit) depending on local sub-domain requirements. An implicit-explicit predictor-corrector scheme, using element-cut partitioning, was developed for Newmark method [12] and its stability analysis was provided using Energy method [14,15]. This approach was later extended to non-linear problems [13] with a proof for convergence [16] and later augmented to control numerical dissipation of unwanted high frequency oscillations using an approach similar to the HHT-$\alpha$ method [17]. Subcycling [18–20] is another approach that use different time-steps in different nodal partitions and is stable as long as the Courant condition is satisfied within respective partitions. These techniques however, cannot handle multiple sub-domains or heterogeneous material properties. Also, the stability/accuracy analysis and implementation of these methods is significantly involved [21]. GC method [22–25] (by Gravouil and Combescure) integrated DDM and implicit-explicit multi-time-stepping enforced by continuity of velocities across interfaces, allowing multiple sub-domains to be integrated independent of each other. It incorporates constituent algorithms from Newmark family and mortar methods to couple distinct discretizations of a continuous domain. It recommends hybrid discretization for the interface of Lagrange multipliers [26], inherently producing several additional unknowns [27]. It was later extended to include Simo, Krenk, HHT-$\alpha$ [28,29] and Multigrid methods [30]. GC method is however unstable with increasing numerical dissipation for large time-step ratios, making it impractical for simulating large scale problems with highly disparate discretizations. PH method [31,32] (by Prakash and Hjelmstad) later introduced a modification to the GC method, making it unconditionally stable for large time-step ratios. It ensures zero interface dissipation for Newmark algorithm, as long as local equilibrium is satisfied, and validates global energy balance via Energy method. PH method is however designed to handle only two sub-domains at a time. Accordingly, it requires a complex and a computationally involved recursive algorithm [32,33] to compute the global solution for multiple sub-domains.

Our goal, through this paper, is to provide a precise of techniques used in multiple grid multiple time-scale (MGMT) coupling, while overcoming some of the limitations mentioned above. That is, we will present a consistent MGMT approach that enables concurrent simulation of multiple sub-domains – with exclusively independent discretizations, energy preserving interfaces and easy to implement solution algorithm. We will provide a comprehensive description of the involved techniques, their derivation and a step-by-step procedure to obtain the global solution for a linear structural dynamic system. We will also provide a rigorous proof of stability using Energy method. Finally, through some numerical examples, we will provide a broad analysis of MGMT results wherein we will evaluate the numerical stability, accuracy and efficiency of the proposed algorithm.

2. Domain Decomposition

Domain decomposition methods (DDM) are among the most efficient and reliable techniques for the solution of engineering applications using FEM. One of the greatest advantages of this method is that the domain under analysis can be decomposed into several component sub-domains, which can then be formulated numerically, modeled and solved independently. Global solution is then obtained by assembling these sub-domains enforced by an interface condition, for example – continuity of a field variable. This allows us to selectively discretize our governing equation, both in space and time, for different sub-domains. In the following section, we present the DDM approach for deriving coupled equations for component sub-domains augmented with an appropriate interface condition. We limit our discussion to structural dynamic systems with damping and linear constitutive relationships.

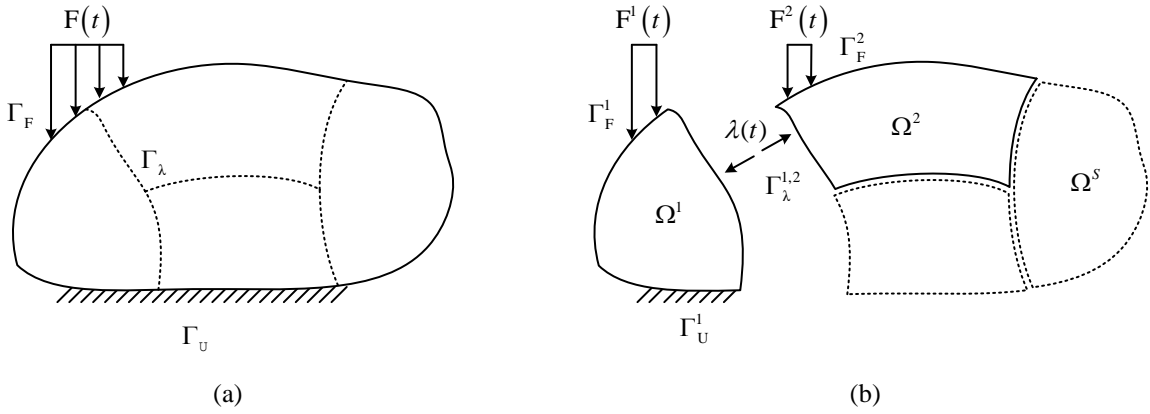

Fig. 1. (a) Structural domain under analysis (b) Decomposed sub-domains with inherited boundary conditions and augmented interface reactions.

Consider a continuous domain $\Omega$ with prescribed displacements over $\Gamma_U$ and prescribed tractions over $\Gamma_F$ as shown in Fig. 1 (a). Using finite element discretization in space, the governing equation / equation of motion can be expressed as:

$$M\ddot{U}(t) + C\dot{U}(t) + KU(t) = F(t) \tag{1}$$

$$\left. \begin{array}{l} U(0) = U_0 \\ \dot{U}(0) = \dot{U}_0 \end{array} \right\} \tag{2}$$

$$\left. \begin{array}{ll} U(t) = U_D(t) & \text{on} \quad \Gamma_U \\ F(t) = F_D(t) & \text{on} \quad \Gamma_F \end{array} \right\} \tag{3}$$

Here, M represents the mass, C represents damping and K represents the stiffness. Primary unknowns are displacements, velocities and accelerations; as represented by U(t), $\dot{U}$(t) and $\ddot{U}$(t) respectively. Initial displacements and velocities are defined by $U_0$ and $\dot{U}_0$ whereas $U_D$(t) and $F_D$(t) represents the prescribed displacement and force boundary conditions. Energy balance equation corresponding to Eq. (1) can be expressed by pre-multiplying $\dot{U}^T$, followed by integration:

$$\frac{d}{dt}\left(\frac{1}{2}\dot{U}^T M \dot{U} + \frac{1}{2} U^T K U\right) = \dot{U}^T F - \dot{U}^T C \dot{U} \tag{4}$$

Using DDM, let us now decompose our original domain into S component sub-domains such that $\Gamma_\lambda^{i,j}$ represents the dividing interface between $\Omega^i$ and $\Omega^j$, as shown in Fig. 1 (b). If Lagrange multipliers in the form of fluxes or tractions are used to represent interface reactions between connecting sub-domains the resulting equation of motion is obtained as:

$$M^i \ddot{U}^i(t) + C^i \dot{U}^i(t) + K^i U^i(t) = F^i(t) - {L^i}^T \lambda(t) \quad (i = 1, 2 \ldots S) \tag{5}$$

And the corresponding energy balance equation is expressed as:

$$\frac{d}{dt}\left(\frac{1}{2}\dot{U}^{i^T} M^i \dot{U}^i + \frac{1}{2} U^{i^T} K^i U^i\right) = \dot{U}^{i^T} F^i - \dot{U}^{i^T} C^i \dot{U}^i - \dot{U}^{i^T} {L^i}^T \lambda \quad (i = 1, 2 \ldots S) \tag{6}$$

Continuity of an unknown variable (say x) across sub-domain interface $\Gamma_\lambda^{i,j}$ can be represented by a linear constraint equation such as:

$$L^i \, x^i(t) = L^j \, x^j(t) \quad or \quad L^i \, x^i(t) - L^j \, x^j(t) = 0 \tag{7}$$

In the above expressions, interface connectivity L represents a Boolean or a multi-constraint operator that is respectively zero and non-zero for interior and interface degrees of freedom. It is used to project required constraints over the selected set of interface degrees of freedom. Variable x may represent displacement (d-continuity), velocity (v-continuity) or acceleration (a-continuity). Since Eqs. (5) and (7) are coupled, it is necessary to ensure that enforcing continuity of variable x does not influence the global energy balance. It is also important to ensure that the time integration of Eq. (5), augmented with the interface condition, yields a stable solution within component sub-domains and globally. These requirements will help us select an appropriate variable for the interface condition. Accordingly, we enforce the following constraints on Eq. (6):

a.  Energy is conserved within respective sub-domains, that is local equilibrium is satisfied.
b.  Interface energy, produced as a result of interface reactions, is equal to zero.

These constraints are reasonable since the global solution (obtained as an assembly of local solutions) can be stable if and only if the solution is stable within respective sub-domains. Also, if the interface energy contributions are positive, numerical integration of Eq. (5) will eventually explode with an unstable solution and if the interface energy is negative, artificial damping will be introduced across sub-domain interfaces. In MGMT method, DD is only a virtual decomposition of a continuous domain; therefore the dividing interfaces do not represent a physical feature within the system. Accordingly, any energy produced as a result of introducing Lagrange multipliers must be zero. This will also ensure stable gluing of adjacent sub-domains and seamless propagation of structural waves across connecting interfaces. Global energy balance equation obtained from component sub-domain contributions can be expressed as:

$$\sum_{i=1}^{S} \frac{d}{dt}\left( \frac{1}{2} \dot{U}^{i^T} M^i \dot{U}^i + \frac{1}{2} U^{i^T} K^i U^i \right) = \sum_{i=1}^{S} \left( \dot{U}^{i^T} F^i - \dot{U}^{i^T} C^i \dot{U}^i - \dot{U}^{i^T} L^{i^T} \lambda \right) \tag{8}$$

By enforcing the constraints discussed above, we obtain the following condition:

$$\sum_{i=1}^{S} \dot{U}^{i^T} L^{i^T} \lambda = 0 \tag{9}$$

Equation (9) represents the necessary condition to ensure global energy balance; specifically it suggests that the interface energy should be identically zero. Comparing this result with Eq. (7), we see that zero interface energy naturally yields continuity of velocities across sub-domain interfaces. The final set of equations for coupled sub-domains can now be represented as:

$$M^i \ddot{U}^i(t) + C^i \dot{U}^i(t) + K^i U^i(t) = F^i(t) - L^{i^T} \lambda(t) \quad (i = 1, 2 \ldots S) \tag{10}$$

$$\sum_{i=1}^{S} L^i \dot{U}^i(t) = 0 \tag{11}$$

### 3. Multiple grid coupling

#### 3.1. The mortar finite element method

One clear advantage of using DDM is the ability to choose independent resolutions in component sub-domains so as to obtain a discretization that is well suited to the local characteristics of the solution to be approximated. Mortar finite element method is an interface discretization technique that is used to couple sub-domain grids across separate spatial resolutions. In this method, node-to-node connectivity may or may not exist at the interface between adjacent sub-domains. Coupling is achieved in the form of point constraints which are enforced by introducing Lagrange multipliers, chosen wisely to preserve the accuracy of the solution, [4,5]. As an example, consider a 2D domain which has been decomposed into two non-overlapping, non-conforming sub-domains ($\Omega^1$ and $\Omega^2$) as shown in Fig. 2 (a). Let them be spatially discretized using 4-node quadrilateral elements with standard bilinear shape functions and let the interface reactions (fluxes or tractions) be represented by Lagrange multipliers ($\lambda$). These unknowns, discretized over the interface, may be approximated using similar bilinear shape functions $N^\lambda$ as shown in Fig. 2 (b).

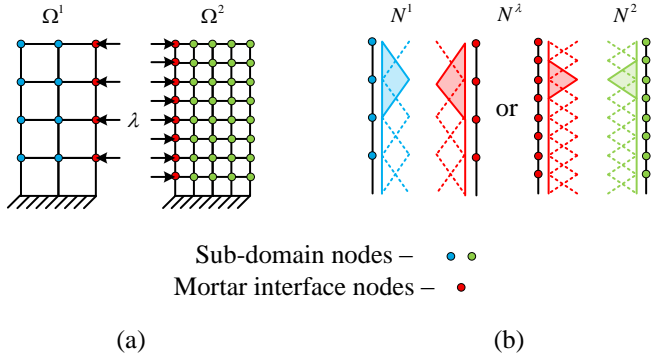

Sub-domain nodes –
Mortar interface nodes –

(a) (b)

Fig. 2. (a) Non-overlapping DD with independent (non-conforming) discretization in component sub-domains (b) Interface elements and corresponding Lagrange multiplier space.

| Region | dof |
|---|---|
| $\Omega^1$ | a |
| $\Omega^2$ | b |
| $\Gamma^1_\lambda$ | m ( < a ) |
| $\Gamma^2_\lambda$ | n ( < b ) |
| $\Gamma^{1,2}_\lambda$ | k |

Table 1. Labels for sub-domain and interface degree of freedom (dof).

Here, if we choose coarse grid discretization to represent the Lagrange multipliers, the adjacent interface on $\Omega^1$ is referred to as the non-mortar interface, whereas the interface on $\Omega^2$ is referred to as the mortar (glued) interface. If multiple sub-domains exist at a node on a mortar element, the usual (linear interpolation) shape functions associated with this particular node are replaced by a constant part [34]. The choice of mortar element discretization (for unknown Lagrange multipliers) can be obtained from either sub-domain or independently. The total number of dof associated with the mortar interface should not be too rich in space so that they over constrain the coarse grid or too weak so that the constraints are not well enforced [35]. GC method uses a hybrid approach in which the mortar dof k = (m + n) – c, where c represents total number of common/coincident dof. In either case, interface connectivity constraints [34] between adjacent sub-domains are obtained as:

$$P^1_{km} = \int_{\Gamma_\lambda} N^\lambda_k N^1_m \qquad (12)$$

$$P^2_{kn} = \int_{\Gamma_\lambda} N^\lambda_k N^2_n \qquad (13)$$

Hence, the central idea behind mortar element methods is to decompose the domain of our interest into non-overlapping sub-domains (using node-cut partitioning) and impose a weak continuity condition across the interface by requiring that the jump of the solution is orthogonal to a suitable Lagrange multiplier space [35,36]. In the analysis of large structures, this approach has two noteworthy advantages. First, the discretization of domain can be selectively improved in localized regions, such as around corners or other features where error in solution is likely to be greatest. This will allow for greater accuracy without the computational burden associated with improving the discretization over the entire global domain. Another practical benefit of this method is that it can be utilized to connect independently modeled and analyzed sub-structures in a large problem. For example, in analysis of an automobile the external framework and the chassis may be modeled independently by different engineers. It is unlikely that sub-structures like these, modeled independently, will have exact node-connectivity at the interface when assembled as a whole for a complete analysis. Transition meshes can be used in the vicinity of the interface, but this would require re-meshing and it would also make the analysis more complex and expensive. Non-conforming mortar methods completely circumvent this difficulty. Also in comparison with the finite element tearing and interconnecting method (FETI) [6], mortar methods have following advantages [37]:

a. It satisfies the compatibility condition between discrete spaces.
b. It provides an inf-sup condition that is independent of the discretization parameter.
c. It results in algebraic systems with well-conditioned matrices.

### 3.2. Interface connectivity

In this section we will discuss how interface constraints can be implemented using non-overlapping DDM and node-cut partitioning. Consider the decomposition of a continuous region $\Omega$ (2D) into two sub-domains, $\Omega^1$ and $\Omega^2$ which are joined together by introducing Lagrange multipliers (interface reactions) at the dividing interface. Since we are using node-cut partitioning, the interface for Lagrange multipliers is a 1D segment. As discussed earlier, the choice for mortar interface discretization is arbitrary but fixed. In our derivation for multiple grid coupling, we choose the coarse grid interface as the non-mortar surface; therefore the dof assigned to the interface of Lagrange multipliers is equal to the $\Omega^1$ interface dof, that is k = m.

#### 3.2.1. Coupling conforming grids

An interface between adjacent sub-domains is said to be conforming if and only if m = n and the coordinates for corresponding dof are coincident.

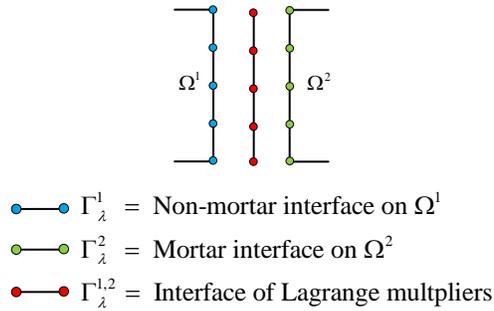

○──● $\Gamma^1_\lambda$ = Non-mortar interface on $\Omega^1$
○──● $\Gamma^2_\lambda$ = Mortar interface on $\Omega^2$
●──● $\Gamma^{1,2}_\lambda$ = Interface of Lagrange multpliers

Fig. 3. Conforming sub-domains (k = m = n).

In this particular case, the interface connectivity operator L is simply a Boolean projection matrix with 1's and 0's assigned to sub-domain interface dof and interior dof respectively. Let $B^1$ and $B^2$ represent the Boolean projection matrix for $\Omega^1$ and $\Omega^1$ respectively. Then the sub-domain dof, located on interface $\Gamma^i_\lambda$, can be selectively projected as:

$$x^1_m\big|_{\Gamma^1_\lambda} = B^1_{ma} x^1_a \qquad \therefore L^1 = B^1 \tag{14}$$

$$x^2_n\big|_{\Gamma^2_\lambda} = B^2_{nb} x^1_b \qquad \therefore L^2 = B^2 \tag{15}$$

In the above expressions, x represents the nodal dof, for example in Eq. (11), x represents the sub-domain velocity vector. As an example for implementing and constructing the Boolean projection matrix, consider the illustration shown in Fig. 4.

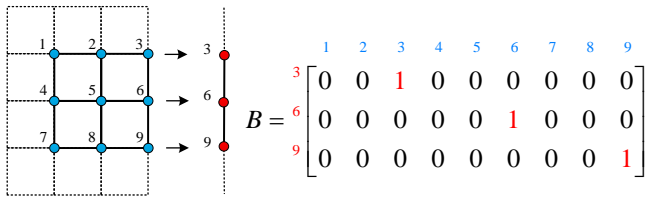

Fig. 4. Reference sub-domain and corresponding Boolean projection matrix.

*3.2.2. Coupling non-conforming grids*

An interface between adjacent sub-domains is said to be non-conforming if m ≠ n and/or the coordinates of corresponding dof are not coincident. In this case, the total number of Lagrange multipliers is selected as k = m, as shown in Fig. 5.

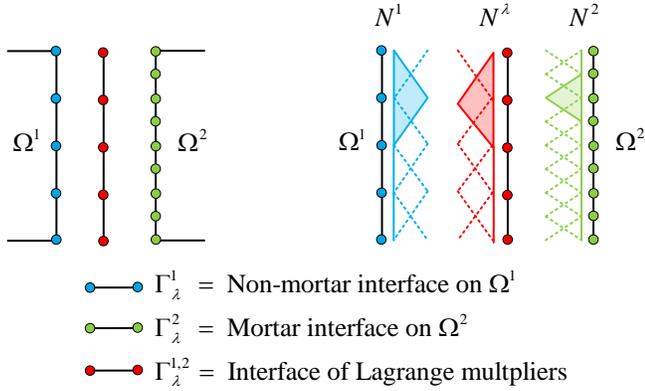

Fig. 5. Non-conforming sub-domains.

Since the component sub-domains are discretized using 4-node (bilinear) quadrilateral elements, we assume equivalent shape functions over the interface of Lagrange multipliers. In order to connect non-conforming grids, first a set of interface dof is projected using the Boolean projection matrix, and then the interface constraints are computed using Eq. (12) or (13). Accordingly the interface connectivity is expressed as:

$$L^1_{ka} = P^1_{km} B^1_{ma} \tag{16}$$

$$L^2_{kb} = P^2_{kn} B^2_{nb} \tag{17}$$

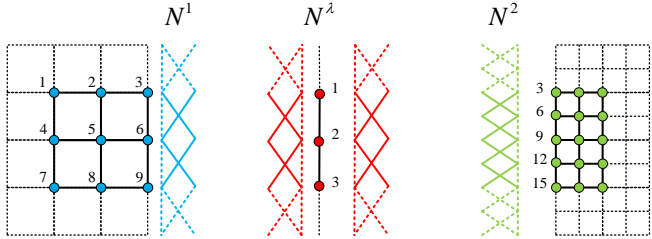

Fig. 6. Reference sub-domains, Interface of Lagrange multipliers and corresponding shape functions.

As an example for implementing and constructing interface connectivity between non-conforming grids, consider the illustration shown in Fig. 6. In this particular example, we will assume that the element edge from $\Omega^1$ and $\Omega^2$ is of 1 and 0.5 unit length respectively. Shape functions used across either interfaces are assumed to be bilinear functions (standard hat functions). The computation for product of shape functions across the interface length, as described in Eqs. (12) and (13), is performed using the Trapezoidal rule [37]. However, integration may also be performed over each element edge using the Quadrature rule [34]. In either case, multiple constrains used in interface connectivity are obtained as:

$$P^1 = \begin{bmatrix} \int N^\lambda_1 N^1_3 & \int N^\lambda_1 N^1_6 & \int N^\lambda_1 N^1_9 \\ \int N^\lambda_2 N^1_3 & \int N^\lambda_2 N^1_6 & \int N^\lambda_2 N^1_9 \\ \int N^\lambda_3 N^1_3 & \int N^\lambda_3 N^1_6 & \int N^\lambda_3 N^1_9 \end{bmatrix} = \begin{bmatrix} 0.33335 & 0.16665 & 0 \\ 0.16665 & 0.66665 & 0.16665 \\ 0 & 0.16665 & 0.33335 \end{bmatrix} \tag{18}$$

$$P^2 = \begin{matrix} 1 \\ 2 \\ 3 \end{matrix} \begin{bmatrix} \int N_1^\lambda N_3^2 & \cdots & \cdots & \cdots & \int N_1^\lambda N_{15}^2 \\ \int N_2^\lambda N_3^2 & \cdots & \cdots & \cdots & \int N_2^\lambda N_{15}^2 \\ \int N_3^\lambda N_3^2 & \cdots & \cdots & \cdots & \int N_3^\lambda N_{15}^2 \end{bmatrix} = \begin{matrix} 1 \\ 2 \\ 3 \end{matrix} \begin{bmatrix} 0.20835 & 0.25 & 0.04165 & 0 & 0 \\ 0.04165 & 0.25 & 0.4167 & 0.25 & 0.04165 \\ 0 & 0 & 0.04165 & 0.25 & 0.20835 \end{bmatrix} \quad (19)$$

with column indices 3, 6, 9, 12, 15.

And the corresponding constraint equations are expressed as:

$$\begin{aligned} 0.33335 x_3^1 + 0.16665 x_6^1 &= 0.20835 x_3^2 + 0.25 x_6^2 + 0.04165 x_9^2 \\ 0.16665 x_3^1 + 0.66665 x_6^1 + 0.16665 x_9^1 &= 0.04165 x_3^2 + 0.25 x_6^2 + 0.4167 x_9^2 + 0.25 x_{12}^2 + 0.04165 x_{15}^2 \\ 0.16665 x_6^1 + 0.33335 x_9^1 &= 0.04165 x_9^2 + 0.25 x_{12}^2 + 0.20835 x_{15}^2 \end{aligned} \quad (20)$$

Hence multiple grid coupling between conforming and non-conforming sub-domains can be achieved using Eqs. (14) - (15) and Eqs. (16) - (17). Note that the interface continuity condition for coupled sub-domains, Eqs. (7) and (11), is expressed as the sum of component sub-domain contributions. Accordingly, the interface connectivity matrix (L) should be assigned a positive bias (+L) for non-mortar dof and a negative bias (-L) for the mortar dof.

### 4. Multiple time-scale coupling

#### 4.1. Newmark time integration

Once discretized in space, the evolution of unknown quantities is obtained by marching in time; that is by direct integration of governing equations. This is when space discretized equations are further discretized in time domain. In this section we will briefly discuss the Newmark time integration method (implicit and explicit) and its implementation algorithm as suited for MGMT coupling. Equation (1) is discretized (temporally) such that the time interval of interest is divided into N sub-steps of size $\Delta t$ and the equilibrium equation is enforced at discrete instants of time $t_n$ where $n \in \{0, 1 \ldots N-1\}$ such that:

$$M\ddot{U}_{n+1} + C\dot{U}_{n+1} + KU_{n+1} = F_{n+1} \quad (21)$$

Using Newmark method, displacement and velocity updates are expressed as:

$$\left. \begin{aligned} U_{n+1} &= U_{n+1}^{pre} + \Delta t^2 \beta \ddot{U}_{n+1} \\ U_{n+1} &= U_n + \Delta t \dot{U}_n + \Delta t^2 (0.5 - \beta) \ddot{U}_n + \Delta t^2 \beta U_{n+1} \end{aligned} \right\} \quad (22)$$

$$\left. \begin{aligned} \dot{U}_{n+1} &= \dot{U}_{n+1}^{pre} + \Delta t \gamma \ddot{U}_{n+1} \\ \dot{U}_{n+1} &= \dot{U}_n + \Delta t (1-\gamma) \ddot{U}_n + \Delta t \gamma \ddot{U}_{n+1} \end{aligned} \right\} \quad (23)$$

Substituting these equations back into Eq. (21), the acceleration-form of Newmark method may be obtained as:

$$\mathbf{K} \ddot{U}_{n+1} = \mathbf{F_{n+1}} \quad (24)$$

Where:
The effective stiffness matrix $\mathbf{K} = \{M + \Delta t \gamma C + \Delta t^2 \beta K\}$  (25)
And the effective load vector:

$$\begin{aligned} \mathbf{F_{n+1}} &= F_{n+1} - C\left\{\dot{U}_n + \Delta t(1-\gamma)\ddot{U}_n\right\} \\ &\quad - K\left\{U_n + \Delta t\dot{U}_n + \Delta t^2(0.5-\beta)\ddot{U}_n\right\} \\ \therefore \mathbf{F_{n+1}} &= F_{n+1} - C\dot{U}_{n+1}^{pre} - KU_{n+1}^{pre} \end{aligned} \quad (26)$$

Equation (24) can be solved for unknown accelerations by factorizing $\mathbf{K} = \mathbf{LDL^T}$ (only once for linear systems with a uniform time-step). Displacements and velocities can then be obtained from Eqs. (22) and (23) respectively. A pseudo-code to advance the solution of an input system of equations from time-step n ∈ {0, 1 …N-1} to (n + 1) using the Newmark method can now be expressed by Algorithm 1.

| Input | |
|---|---|
| M   C   K | $F_n$ |
| $\ddot{U}_n$   $\dot{U}_n$   $U_n$ | $\Delta t$   $\beta$   $\gamma$ |
| Black-box 1: Newmark integration | |
| 1) Compute effective stiffness matrix using Eq. (25) <br> 2) Compute effective load vector using Eq. (26) <br> 3) Solve for updated accelerations using Eq. (24) <br> 4) Compute updated displacements using Eq. (22) <br> 5) Compute updated velocities using Eq. (23) | |
| Output | |
| $\ddot{U}_{n+1}$   $\dot{U}_{n+1}$   $U_{n+1}$ | |

Algorithm 1: Time integration using Newmark method – acceleration form (Black-box 1).

The Newmark algorithm (for an undamped case) is stable if γ ≥ 1/2 and is unconditionally stable if β ≥ 1/4(γ + 1/2)². With appropriate expression for β and γ, the Newmark method is further classified as [38]:

a.  Implicit constant average acceleration method (also known as Trapezoidal rule) for γ = 1/2 and β = 1/4
b.  Implicit linear acceleration method (corresponds to θ = 1 in the Wilson method) for γ = 1/2 and β = 1/6
c.  Explicit central difference method for γ = 1/2 and β = 0

Equations (21), (22) and (23) can be expressed together in a compact equation as follows:

$$\tilde{M}\tilde{U}_{n+1} = \tilde{F}_{n+1} - \tilde{N}\tilde{U}_n \quad (27)$$

Where:

$$\tilde{M} = \begin{bmatrix} M & C & K \\ -\Delta t\gamma & I & 0 \\ -\Delta t^2\beta & 0 & I \end{bmatrix} \quad (28)$$

$$\tilde{N} = \begin{bmatrix} M & C & K \\ -\Delta t(1-\gamma)I & -I & 0 \\ -\Delta t^2(0.5-\beta) & -\Delta t I & -I \end{bmatrix} \quad (29)$$

$$\tilde{F}_{n+1} = \begin{bmatrix} F_{n+1} & 0 & 0 \end{bmatrix}^T \quad (30)$$

$$\tilde{U}_n = \begin{bmatrix} \ddot{U}_n & \dot{U}_n & U_n \end{bmatrix}^T \tag{31}$$

### 4.2. Interface condensation

In this section we will derive the fully discretized equations of motion for decomposed sub-domains, coupled over multiple time scales. We shall use the DD from Fig. 1 (b). It is assumed that every component sub-domain is discretized (spatially and temporally) independent of each other and multiple grid coupling between adjacent sub-domains is achieved using mortar finite elements. Let the global time-step for evolving MGMT equations from time $t_n$ where $n \in \{0,1,\ldots N-1\}$ be $\Delta T$. For reference purpose, we also assume that $\Omega^1$ is integrated at the global time-step $\Delta T$. Every other component sub-domain can be integrated with time-steps $\Delta t^i$ such that $(T + \Delta T) = (T + \xi^i \Delta t^i)$. In order to obtain global solution at synchronous time instants it is necessary to ensure that the time-step ratio $\xi^i$ is an integer factor of the global time-step $\Delta T$. Accordingly, we will use the following descriptions:

| Sub-domain | Time-step ($\Delta t$) | Algorithmic parameters | Time-step ratio ($\xi = \Delta T / \Delta t$) | Intermediate step counter ($\eta$) |
|---|---|---|---|---|
| $\Omega^1$ | $\Delta t^1 = \Delta T$ | $\beta^1, \gamma^1$ | $\xi^1 = 1$ | $\eta^1 = 0, 1$ |
| $\Omega^2$ | $\Delta t^2$ | $\beta^2, \gamma^2$ | $\xi^2$ | $\eta^2 = 0, 1, 2 \ldots \xi^2$ |
| $\vdots$ | $\vdots$ | $\vdots$ | $\vdots$ | $\vdots$ |
| $\Omega^S$ | $\Delta t^S$ | $\beta^S, \gamma^S$ | $\xi^S$ | $\eta^S = 0, 1, 2, 3 \ldots \xi^S$ |

Table 2. Sub-domain time-stepping parameters.

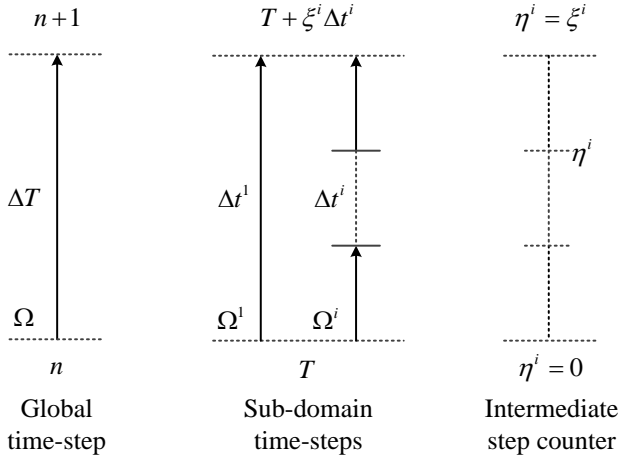

Fig. 7. MGMT time-stepping.

Using Newmark method, fully discretized equations of motion for component sub-domains along with the interface condition enforced at $\eta^i = \xi^i$ (global time-steps) are expressed as:

$$M^i \ddot{U}^i_{\eta^i} + C^i \dot{U}^i_{\eta^i} + K^i U^i_{\eta^i} + L^{i^T} \lambda_{\eta^i} = F^i_{\eta^i} \quad (i = 1, 2 \ldots S) \tag{32}$$

$$\sum_{i=1}^{S} L^i \dot{U}^i_{\xi^i} = 0 \tag{33}$$

In Eq. (32), unknown interface reactions (Lagrange multipliers) are computed at every intermediate time-step, for every component sub-domain. This is computationally very expensive, especially in cases where time-step ratio between sub-domains is large. Intermediate interface reactions (at $0 < \eta^i < \xi^i$) can however be condensed, so that Eq. (32) can be expressed in terms of Lagrange multipliers at global time-steps, that is in terms of $\lambda_n$ ($\eta^i = 0$) and $\lambda_{n+1}$ ($\eta^i = \xi^i$). Interface condensation is performed using an energy preserving approach that is similar to PH method [31]. This approach makes the interface computation $\xi^i$ times faster than GC method [23], making it computationally superior than the GC method. The drawback in using PH method is that it allows coupling only two sub-domains at a time. Hence for multiple sub-domain coupling, one needs to implement a recursive solution algorithm [32] that can solve component sub-domains one pair at a time. This again restricts the computational efficiency of multiple time-scale coupling. Also, in PH method, the Lagrange multipliers are expressed in terms of 'unbalanced' interface reactions, which are obtained from the coarse time-step sub-domain. This makes the fine time-step sub-domain dependent on the other, making it impossible to solve multiple sub-domains concurrently (using parallel processing). We however use the energy preserving approach with a slight modification such that the intermediate interface reactions are expressed only in the terms of known ($\lambda_n$) and unknown ($\lambda_{n+1}$) Lagrange multipliers. We will also show that our modification results in exclusively independent sub-domains, so that they may be solved concurrently (using parallel processing), further improving computational efficiency. We begin with kinematic decomposition, Eq. (34), which will enable us to define the equation of motion independently under the action of external forces and under the response of interface reactions.

$$\tilde{U} = \tilde{V} + \tilde{W} \tag{34}$$

Where the notation (~) was defined earlier in Eq. (31). Using Eq. (34), equilibrium of $\Omega^1$ can be expressed as:

$$\left(M^1 \ddot{V}^1_{n+1} + C^1 \dot{V}^1_{n+1} + K^1 V^1_{n+1}\right) + \left(M^1 \ddot{W}^1_{n+1} + C^1 \dot{W}^1_{n+1} + K^1 W^1_{n+1}\right) = F^1_{n+1} - L^{1^T} \lambda_{n+1} \tag{35}$$

Note for $\Omega^1$, $\xi^1 = 1$ which corresponds to instant (n + 1). Equation (35) can now be decomposed into two equations as follows:

$$M^1 \ddot{V}^1_{n+1} + C^1 \dot{V}^1_{n+1} + K^1 V^1_{n+1} = F^1_{n+1} \tag{36}$$

$$M^1 \ddot{W}^1_{n+1} + C^1 \dot{W}^1_{n+1} + K^1 W^1_{n+1} = -L^{1^T} \lambda_{n+1} \tag{37}$$

Equation (36) represents the equilibrium of $\Omega^1$ and its contribution to the kinematic quantities under the action of external forces only. Equation (37) on other hand represents the equilibrium under the action of Lagrange multipliers or unknown interface reactions. In order to avoid interface dissipation across intermediate time-steps, the equilibrium of $\Omega^1$ is also enforced at every $\eta^i$ instant by requiring that the combined residual of Eqs. (36) and (37) is equal to zero. That is:

$$\left(\delta f\right)^1_{\eta^i} + \left(\delta l\right)^1_{\eta^i} = 0 \quad (i = 1, 2 \ldots S) \tag{38}$$

For i = 1, $\eta$ = 1, Eq. (38) is identical to solving the equilibrium of $\Omega^1$ as defined by Eq. (32). However, for $1 < i \leq S$ and $\eta > 1$, it is equivalent to enforcing the equilibrium of $\Omega^1$ across intermediate time-steps. The first term in Eq. (38) represents the free residuals obtained under the action of external forces and the second term represents the link residuals obtained under the action of unknown interface reactions:

$$\therefore \left(\delta f\right)^1_{\eta^i} = F^1_{\eta^i} - M^1 \ddot{V}^1_{\eta^i} - C^1 \dot{V}^1_{\eta^i} - K^1 V^1_{\eta^i} \quad (i = 2 \ldots S) \tag{39}$$

$$\therefore \left(\delta l\right)^1_{\eta^i} = -L^{1^T} \lambda_{\eta^i} - M^1 \ddot{W}^1_{\eta^i} - C^1 \dot{W}^1_{\eta^i} - K^1 W^1_{\eta^i} \quad (i = 2 \ldots S) \tag{40}$$

Before enforcing Eq. (38) let us define necessary interpolation functions for obtaining intermediate quantities used in Eq. (39) and (40):

$$\tilde{V}^1_{\eta^i} = \left(1 - \frac{\eta^i}{\xi^i}\right) \tilde{U}^1_n + \left(\frac{\eta^i}{\xi^i}\right) \tilde{V}^1_{n+1} \tag{41}$$

$$\tilde{W}^1_{\eta^i} = \left(\frac{\eta^i}{\xi^i}\right) \tilde{W}^1_{n+1} \tag{42}$$

$$F^1_{\eta^i} = \left(1 - \frac{\eta^i}{\xi^i}\right) F^1_n + \left(\frac{\eta^i}{\xi^i}\right) F^1_{n+1} \tag{43}$$

Equation (36) can be solved independently using the Newmark method, Algorithm 1, yielding the solution vector at time instant (n + 1) or ($\eta^1 = \xi^1 = 1$). Equation (39) can now be expressed as:

$$\therefore (\delta f)^1_{\eta^i} = \left\{\left(1 - \frac{\eta^i}{\xi^i}\right) F^1_n + \left(\frac{\eta^i}{\xi^i}\right) F^1_{n+1}\right\} - M^1 \left\{\left(1 - \frac{\eta^i}{\xi^i}\right) \ddot{U}^1_n + \left(\frac{\eta^i}{\xi^i}\right) \ddot{V}^1_{n+1}\right\}$$
$$- C^1 \left\{\left(1 - \frac{\eta^i}{\xi^i}\right) \dot{U}^1_n + \left(\frac{\eta^i}{\xi^i}\right) \dot{V}^1_{n+1}\right\} - K^1 \left\{\left(1 - \frac{\eta^i}{\xi^i}\right) U^1_n + \left(\frac{\eta^i}{\xi^i}\right) V^1_{n+1}\right\} \tag{44}$$

$$\therefore (\delta f)^1_{\eta^i} = \left(1 - \frac{\eta^i}{\xi^i}\right) \left\{F^1_n - M^1 \ddot{U}^1_n - C^1 \dot{U}^1_n - K^1 U^1_n\right\}$$
$$+ \left(\frac{\eta^i}{\xi^i}\right) \left\{F^1_{n+1} - M^1 \ddot{V}^1_{n+1} - C^1 \dot{V}^1_{n+1} - K^1 V^1_{n+1}\right\} \tag{45}$$

Using Eqs. (32) and (36) we have:

$$(\delta f)^1_{\eta^i} = \left(1 - \frac{\eta^i}{\xi^i}\right) L^{1^T} \lambda_n \tag{46}$$

Using Eq. (42), link residuals from Eq. (40) can be expressed as:

$$(\delta l)^1_{\eta^i} = -L^{1^T} \lambda_{\eta^i} - M^1 \left(\frac{\eta^i}{\xi^i}\right) \ddot{W}^1_{n+1} - C^1 \left(\frac{\eta^i}{\xi^i}\right) \dot{W}^1_{n+1} - K^1 \left(\frac{\eta^i}{\xi^i}\right) W^1_{n+1}$$
$$\therefore (\delta l)^1_{\eta^i} = -L^{1^T} \lambda_{\eta^i} - \left(\frac{\eta^i}{\xi^i}\right) \left\{M^1 \ddot{W}^1_{n+1} + C^1 \dot{W}^1_{n+1} + K^1 W^1_{n+1}\right\} \tag{47}$$

Further using Eq. (37):

$$(\delta l)^1_{\eta^i} = -L^{1^T} \lambda_{\eta^i} + \left(\frac{\eta^i}{\xi^i}\right) L^{1^T} \lambda_{n+1} \tag{48}$$

Enforcing Eq. (38) we get:

$$\left(1 - \frac{\eta^i}{\xi^i}\right) L^{1^T} \lambda_n - L^{1^T} \lambda_{\eta^i} + \left(\frac{\eta^i}{\xi^i}\right) L^{1^T} \lambda_{n+1} = 0$$
$$\therefore L^{1^T} \lambda_{\eta^i} = \left(1 - \frac{\eta^i}{\xi^i}\right) L^{1^T} \lambda_n + \left(\frac{\eta^i}{\xi^i}\right) L^{1^T} \lambda_{n+1}$$
$$\therefore \lambda_{\eta^i} = \left(1 - \frac{\eta^i}{\xi^i}\right) \lambda_n + \left(\frac{\eta^i}{\xi^i}\right) \lambda_{n+1} \tag{49}$$

Equation (49) represents the interpolation function to compute intermediate Lagrange multipliers (interface reactions) at every $\eta^i$ step of $\Omega^i$. We can now condense the intermediate Lagrange multipliers and express the equilibrium of $\Omega^i$ as follows:

$$M^i \ddot{U}^i_{\eta^i} + C^i \dot{U}^i_{\eta^i} + K^i U^i_{\eta^i} + \left(\frac{\eta^i}{\xi^i}\right) L^{i^T} \lambda_{n+1} = F^i_{\eta^i} - \left(1 - \frac{\eta^i}{\xi^i}\right) L^{i^T} \lambda_n \tag{50}$$

Comparing Eq. (50) with Eq. (41) from [31], we clearly see that our interface condensation approach does not impose any sub-domain dependency. Intermediate interface reactions are now defined only in terms of unknown Lagrange multipliers at (n+1). This further enables us to solve component sub-domains independent of each other followed by a direct solution of $\lambda_{n+1}$. Using compact notations discussed earlier in Eqs. (27) - (31), global system of fully discretized equations along with displacement and velocity updates can now be expressed as:

$$\tilde{M}^i \tilde{U}^i_{\eta^i} + \left(\frac{\eta^i}{\xi^i}\right) \tilde{L}^i \lambda_{n+1} = \tilde{F}^i_{\eta^i} - \tilde{N}^i \tilde{U}^i_{\eta^i-1} - \left(1 - \frac{\eta^i}{\xi^i}\right) L^{i^T} \lambda_n \quad (i = 1, 2 \ldots S) \tag{51}$$

$$\sum_{i=1}^{S} \tilde{B}^i \tilde{U}^i_{n+1} = 0 \tag{52}$$

Where: $\tilde{B}^i = \begin{bmatrix} 0 & L^i & 0 \end{bmatrix}$

### 5. Stability analysis using energy method

When deriving coupled MGMT equations for component sub-domains, we assumed that the solution is stable if energy is conserved within respective sub-domains (local equilibrium is satisfied) and interface energy, produced as a result of introducing Lagrange multipliers, is equal to zero. This allowed us to define an appropriate variable for enforcing continuity across dividing interfaces. In this section we will show, using Energy method, that the change in energy due to MGMT coupling over time-step ($t_n - t_{n+1}$) or (T – T + ΔT) is identically zero if the time integration scheme is stable within respective sub-domains. We will also show that enforcing continuity of velocities conserves global energy by yielding zero interface energy contributions.

$$[x_n] = x_{n+1} - x_n \qquad \langle x_n \rangle = \frac{x_{n+1} + x_n}{2} \tag{53}$$

Using the undivided forward difference and mean value operators from Eq. (53), the energy balance equation for a single and a continuous domain, identity (9.4.47) from [15], is expressed as:

$$\frac{1}{2}\ddot{U}^T_{n+1} A \ddot{U}_{n+1} - \frac{1}{2}\ddot{U}^T_n A \ddot{U}_n + \frac{1}{2}\dot{U}^T_{n+1} K \dot{U}_{n+1} - \frac{1}{2}\dot{U}^T_n K \dot{U}_n =$$
$$-\left(\gamma - \frac{1}{2}\right) [\ddot{U}_n]^T B [\ddot{U}_n] - \Delta t \langle \ddot{U}_n \rangle^T C \langle \ddot{U}_n \rangle \tag{54}$$

Where:

$$\left. \begin{array}{l} B = M + \Delta t \left(\gamma - \frac{1}{2}\right) C + \Delta t^2 \left(\beta - \frac{\gamma}{2}\right) K \\ A = B + \Delta t \left(\gamma - \frac{1}{2}\right) C \end{array} \right\} \tag{55}$$

Equation (54) for a component sub-domain along with interface reactions can be obtained as:

$$\frac{1}{2}\ddot{U}_{n+1}^T A \ddot{U}_{n+1} - \frac{1}{2}\ddot{U}_n^T A \ddot{U}_n + \frac{1}{2}\dot{U}_{n+1}^T K \dot{U}_{n+1} - \frac{1}{2}\dot{U}_n^T K \dot{U}_n =$$
$$-\left(\gamma - \frac{1}{2}\right)[\ddot{U}_n]^T B [\ddot{U}_n] - \Delta t \langle \ddot{U}_n \rangle^T C \langle \ddot{U}_n \rangle - \frac{1}{\Delta t}[\dot{U}_n]^T L^T [\lambda_n] \quad (56)$$

The system of MGMT equations can now be considered numerically stable if the total energy change (the RHS of Eq. (56)) under no external loads is less than or equal to zero. Hence, if we can show that the interface contributions are less than or equal to zero, the burden of stability relies only upon the stability of component sub-domains, as integrated using Newmark method. Let the interface contributions from Eq. (56) be represented as:

$$\hat{E}_{\Gamma_\lambda} = \frac{1}{\Delta t}[\dot{U}_n]^T L^T [\lambda_n] \quad (57)$$

We shall now assume the DD of a continuous region $\Omega$ into S component sub-domains, Fig. 1 (b), with sub-domain time-stepping parameters described in Fig. 7 and Table 2. Using Eq. (57), the combined interface contributions of component sub-domains is expressed as:

$$\hat{E}_{\Gamma_\lambda} = \sum_{i=1}^{S} \sum_{\eta^i=1}^{\xi^i} \hat{E}_{\Gamma_\lambda^i}$$
$$= \sum_{i=1}^{S} \sum_{\eta^i=1}^{\xi^i} \frac{1}{\Delta t^i}[\dot{U}_{\eta^i-1}^i]^T L^{iT}[\lambda_{\eta^i-1}]$$
$$= \sum_{i=1}^{S} \sum_{\eta^i=1}^{\xi^i} \frac{1}{\Delta t^i}(\dot{U}_{\eta^i}^i - \dot{U}_{\eta^i-1}^i)^T L^{iT}(\lambda_{\eta^i} - \lambda_{\eta^i-1}) \quad (58)$$

Using the relationship derived in Eq. (49) the forward difference of interface reactions in terms of $\lambda_n$ and $\lambda_{n+1}$ is expressed as:

$$[\lambda_{\eta-1}] = \lambda_\eta - \lambda_{\eta-1}$$
$$= \left\{\left(1-\frac{\eta}{\xi}\right)\lambda_n + \left(\frac{\eta}{\xi}\right)\lambda_{n+1}\right\} - \left\{\left(1-\frac{\eta-1}{\xi}\right)\lambda_n + \left(\frac{\eta-1}{\xi}\right)\lambda_{n+1}\right\} \quad (59)$$
$$= \frac{1}{\xi}(\lambda_{n+1} - \lambda_n)$$

Accordingly, Eq. (58) now becomes:

$$\hat{E}_{\Gamma_\lambda} = \sum_{i=1}^{S} \frac{1}{\xi^i \Delta t^i}\left\{\sum_{\eta^i=1}^{\xi^i}(\dot{U}_{\eta^i}^i - \dot{U}_{\eta^i-1}^i)^T L^{iT}(\lambda_{n+1} - \lambda_n)\right\}$$
$$= \sum_{i=1}^{S}(\lambda_{n+1} - \lambda_n)^T \left\{\frac{L^i}{\xi^i \Delta t^i}\sum_{\eta^i=1}^{\xi^i}(\dot{U}_{\eta^i}^i - \dot{U}_{\eta^i-1}^i)\right\} \quad (60)$$

Note that:

$$\sum_{\eta^i=1}^{\xi^i}(\dot{U}_{\eta^i}^i - \dot{U}_{\eta^i-1}^i) = (\dot{U}_1^i - \dot{U}_0^i) + (\dot{U}_2^i - \dot{U}_1^i) + \ldots + (\dot{U}_{\xi^i}^i - \dot{U}_{\xi^i-1}^i)$$
$$= \dot{U}_{\xi^i}^i - \dot{U}_0^i \quad (61)$$
$$= \dot{U}_{n+1}^i - \dot{U}_n^i$$

$$\xi^i \Delta t^i = \Delta T \tag{62}$$

Therefore we have:

$$\begin{aligned}
\hat{E}_{\Gamma_\lambda} &= \sum_{i=1}^{S} (\lambda_{n+1} - \lambda_n)^T \left\{ \frac{L^i}{\Delta T} (\dot{U}_{n+1}^i - \dot{U}_n^i) \right\} \\
&= \frac{(\lambda_{n+1} - \lambda_n)^T}{\Delta T} \sum_{i=1}^{S} L^i (\dot{U}_{n+1}^i - \dot{U}_n^i) \\
&= \frac{(\lambda_{n+1} - \lambda_n)^T}{\Delta T} \sum_{i=1}^{S} L^i \dot{U}_{n+1}^i - L^i \dot{U}_n^i
\end{aligned} \tag{63}$$

Since we have enforced continuity of velocities across sub-domain interfaces, Eq. (11), the above expression yields:

$$\hat{E}_{\Gamma_\lambda} = \frac{(\lambda_{n+1} - \lambda_n)^T}{\Delta T} \left\{ \underbrace{\sum_{i=1}^{S} L^i \dot{U}_{n+1}^i}_{=0} - \underbrace{\sum_{i=1}^{S} L^i \dot{U}_n^i}_{=0} \right\} \tag{64}$$

$$\therefore \hat{E}_{\Gamma_\lambda} = 0 \tag{65}$$

Equation (65) proves that introducing Lagrange multipliers in conjugation with the continuity of velocity constraint, results in identically ZERO energy contributions. Accordingly, the stability of MGMT coupling only depends on the stability of Newmark method in integrating component sub-domains. Hence, as long as the stability requirements are satisfied within the time integration of respective sub-domains, MGMT coupling is stable and energy preserving. Identical stability analysis may be performed by enforcing continuity of displacements or accelerations. Interface energy in these cases is zero under the following conditions [39]:

  a. $2\beta^i - \gamma^i = 0$ (for i = 1,2…S) for displacement continuity (d-continuity)
  b. $\gamma^i$ = constant (for i = 1,2…S) for acceleration continuity (a-continuity)

Accordingly, d or a-continuity can be conveniently modeled to yield zero interface energy, however at the cost of imposing further restrictions on choosing exclusively independent time-integration algorithms for distinct sub-domains.

## 6. Solution algorithm and its FE implementation

In the following section we will describe the implementation, a step-by-step solution procedure, for obtaining the solution of coupled MGMT equations as described by Eqs. (51) and (52). Let us use the DD from Fig. 1 (b) and the time-stepping parameters from Fig. 7 and Table 2. In order to communicate information across component sub-domains at the global time-step $\Delta T$ or ($\eta^i = \xi^i$) instants, equilibrium of $\Omega^i$ (where $\xi^i > 1$) needs to be advanced from $\eta^i = 1$ to $\xi^i$. These sub-steps combined together in a matrix format can be expressed as:

$$\begin{bmatrix} \tilde{M}^i & 0 & 0 & 0 & \vdots & (1/\xi^i)\tilde{L}^i \\ \tilde{N}^i & \tilde{M}^i & 0 & 0 & \vdots & (2/\xi^i)\tilde{L}^i \\ 0 & \ddots & \ddots & 0 & \vdots & \vdots \\ 0 & 0 & \tilde{N}^i & \tilde{M}^i & \vdots & \tilde{L}^i \end{bmatrix} \begin{bmatrix} \tilde{U}_1^i \\ \tilde{U}_2^i \\ \vdots \\ \tilde{U}_{\xi^i}^i \\ \hdashline \lambda_{n+1} \end{bmatrix} = \begin{bmatrix} \tilde{F}_1^i - \tilde{N}^i \tilde{U}_n^i - \left(1 - \frac{1}{\xi^i}\right) L^{i^T} \lambda_n \\ \tilde{F}_2^i - \left(1 - \frac{2}{\xi^i}\right) L^{i^T} \lambda_n \\ \vdots \\ \tilde{F}_{\xi^i}^i \end{bmatrix} \tag{66}$$

Further simplified and expressed in terms of synchronized instants:

$$\tilde{M}^i \tilde{U}^i_{\xi^i} + \tilde{L}^i \lambda_{n+1} = \tilde{F}^i_{\xi^i} \tag{67}$$

Where:

$$\tilde{M}^i = \begin{bmatrix} \tilde{M}^i & 0 & 0 & 0 \\ \tilde{N}^i & \tilde{M}^i & 0 & 0 \\ 0 & \ddots & \ddots & 0 \\ 0 & 0 & \tilde{N}^i & \tilde{M}^i \end{bmatrix} \quad \tilde{U}^i_{\xi^i} = \begin{bmatrix} \tilde{U}^i_1 \\ \tilde{U}^i_2 \\ \vdots \\ \tilde{U}^i_{\xi^i} \end{bmatrix} \tag{68}$$

$$\tilde{L}^i = \begin{bmatrix} (1/\xi^i)\tilde{L}^i \\ (2/\xi^i)\tilde{L}^i \\ \vdots \\ \tilde{L}^i \end{bmatrix} \quad \tilde{F}^i_{\xi^i} = \begin{bmatrix} \tilde{F}^i_1 - \tilde{N}^i \tilde{U}^i_n - \left(1 - \dfrac{1}{\xi^i}\right) L^{i^T} \lambda_n \\ \tilde{F}^i_2 - \left(1 - \dfrac{2}{\xi^i}\right) L^{i^T} \lambda_n \\ \vdots \\ \tilde{F}^i_{\xi^i} \end{bmatrix} \tag{69}$$

For sub-domain $\Omega^1$ (where $\xi^1 = 1$), we have:

$$\tilde{M}^1 \tilde{U}^1_{n+1} + \tilde{L}^1 \lambda_{n+1} = \tilde{F}^1_{n+1} - \tilde{N}^1 \tilde{U}^1_n \tag{70}$$

Final system of equations for coupled sub-domains (discretized independently in space and time) and synchronized at every global time-step $\Delta T$ or $(n + 1)$ can now be expressed as:

$$\left[\begin{array}{cccc|c} \tilde{M}^S & \cdots & 0 & 0 & \tilde{L}^S \\ \vdots & \ddots & \vdots & \vdots & \vdots \\ 0 & \cdots & \tilde{M}^2 & 0 & \tilde{L}^2 \\ 0 & \cdots & 0 & \tilde{M}^1 & \tilde{L}^1 \\ \hline \tilde{B}^S & \cdots & \tilde{B}^2 & \tilde{B}^1 & 0 \end{array}\right] \begin{bmatrix} \tilde{U}^S_{\xi^S} \\ \vdots \\ \tilde{U}^2_{\xi^2} \\ \tilde{U}^1_{n+1} \\ \lambda_{n+1} \end{bmatrix} = \begin{bmatrix} \tilde{F}^S_{\xi^S} \\ \vdots \\ \tilde{F}^2_{\xi^2} \\ \tilde{F}^1_{n+1} - \tilde{N}^1 \tilde{U}^1_n \\ 0 \end{bmatrix} \tag{71}$$

Where: $\tilde{B}^i = \begin{bmatrix} 0 & 0 & \cdots & \tilde{B}^i \end{bmatrix}$

Above set of equations can be conveniently grouped (block wise) as:

$$\left[\begin{array}{c|c} A & b \\ \hline c & d \end{array}\right] \begin{bmatrix} x \\ y \end{bmatrix} = \begin{bmatrix} f \\ g \end{bmatrix} \tag{72}$$

Following the BEC algorithm [40], the LHS coefficient matrix in Eq. (72) can be decomposed using Crout Factorization as follows:

$$\begin{bmatrix} A & b \\ c & d \end{bmatrix} = \begin{bmatrix} A & 0 \\ c & \delta \end{bmatrix} \begin{bmatrix} I & v \\ 0 & 1 \end{bmatrix} \tag{73}$$

Global solution to the original problem at time $(T + \Delta T)$ can now be obtained using the following steps:

### 6.1. Step 1: Solve Av = b

$$\begin{bmatrix} \tilde{\tilde{M}}^S & \cdots & 0 & 0 \\ \vdots & \ddots & \vdots & \vdots \\ 0 & \cdots & \tilde{\tilde{M}}^2 & 0 \\ 0 & \cdots & 0 & \tilde{M}^1 \end{bmatrix} \begin{bmatrix} \tilde{v}^S_{\xi^S} \\ \vdots \\ \tilde{v}^2_{\xi^2} \\ \tilde{v}^1_{n+1} \end{bmatrix} = \begin{bmatrix} \tilde{L}^S \\ \vdots \\ \tilde{L}^2 \\ \tilde{L}^1 \end{bmatrix} \tag{74}$$

Where: $\tilde{v}^i_{\xi^i} = \begin{bmatrix} \tilde{v}^i_1 & \tilde{v}^i_2 & \cdots & \tilde{v}^i_{\xi^i} \end{bmatrix}^T$

Since the LHS (coefficient) matrix and the RHS vector in Eq. (74) remain constant in case of linear algorithms, the above system of equations is computed only once before entering the global time-stepping loop. Starting with $\Omega^1$, where $\xi^1 = 1$ corresponds to (n+1) time instant:

$$\tilde{M}^1 \tilde{v}^1_{n+1} = \tilde{L} \Rightarrow \begin{bmatrix} M^1 & C^1 & K^1 \\ -\Delta T \gamma^1 & I & 0 \\ -\Delta T^2 \beta^1 & 0 & I \end{bmatrix} \begin{bmatrix} \ddot{v}^1_{n+1} \\ \dot{v}^1_{n+1} \\ v^1_{n+1} \end{bmatrix} = \begin{bmatrix} L^T \\ 0 \\ 0 \end{bmatrix} \tag{75}$$

Comparing Eq. (75) with Eq. (27) we see that LHS coefficient matrix represents the Newmark set of equations for $\Omega^1$ with zero initial conditions. The solution vector can hence be obtained by using Black box routine for Newmark method, Algorithm 1. Notice that the RHS (load) vector consists of interface connectivity matrix (Boolean projection matrix in case of conforming sub-domains), the solution array v in Eq. (74) is also referred to as the unit load response matrix. The first index of interface connectivity matrix is equal to the total number of degrees of freedom associated with the interface of Lagrange multipliers, say dof ($\lambda$). Hence, the above system needs to be solved successively by considering one RHS column at a time; that is by loading one interface dof at a time. Accordingly, the solution array will be an aggregate of dof ($\lambda$) columns as follows:

$$\tilde{v}^1_{n+1} = \begin{bmatrix} \overset{j=1}{\downarrow} & & \overset{j=dof(\lambda)}{\downarrow} \\ \ddot{v}^1_{n+1} & \cdots & \ddot{v}^1_{n+1} \\ \dot{v}^1_{n+1} & \cdots & \dot{v}^1_{n+1} \\ v^1_{n+1} & \cdots & v^1_{n+1} \end{bmatrix} \quad \tilde{v}^i_{\eta^i} = \begin{bmatrix} \overset{j=1}{\downarrow} & & \overset{j=dof(\lambda)}{\downarrow} \\ \ddot{v}^i_{\eta^i} & \cdots & \ddot{v}^i_{\eta^i} \\ \dot{v}^i_{\eta^i} & \cdots & \dot{v}^i_{\eta^i} \\ v^i_{\eta^i} & \cdots & v^i_{\eta^i} \end{bmatrix} \tag{76}$$

For component sub-domains where $\xi^i > 1$, the solution is updated after every intermediate time-step until $\xi^i = \eta^i$, that is until the solution is synchronized with (n+1) time instant. Step 1 represented as a sub-domain independent pseudo-code can be illustrated by Algorithm 2.

| for i = 1 to S |
| --- |
| for j = 1 to dof (λ) |
| for k = 1 to $\xi^i$ |
| Input |
| $M^i$ $C^i$ $K^i$    $(k/\xi^i)L^{i^T}[:,j]$ |
| $\ddot{v}_{k-1}$ $\dot{v}_{k-1}$ $v_{k-1}$    $\Delta t^i$ $\beta^i$ $\gamma^i$ |
| Black-box 1: Newmark integration |
| Output |
| $\ddot{v}_k$ $\dot{v}_k$ $v_k \to \tilde{v}_k^i[:,j]$ |

Algorithm 2: Computing unit load response matrix for component sub-domain $\Omega^i$.

### 6.2. Step 2: Compute $\delta = d - cv$

$$\delta = [0] - \begin{bmatrix} \tilde{B}^S & \cdots & \tilde{B}^2 & \tilde{B}^1 \end{bmatrix} \begin{bmatrix} \tilde{v}_{\xi^S}^S \\ \vdots \\ \tilde{v}_{\xi^2}^2 \\ \tilde{v}_{n+1}^1 \end{bmatrix} \tag{77}$$

$$\left.\begin{aligned} \therefore \delta &= -\left[ \tilde{B}^S \tilde{v}_{\xi^S}^S + \ldots + \tilde{B}^2 \tilde{v}_{\xi^2}^2 + \tilde{B}^1 \tilde{v}_{n+1}^1 \right] \\ \therefore \delta &= -\left[ \tilde{B}^S \tilde{v}_{n+1}^S + \ldots + \tilde{B}^2 \tilde{v}_{n+1}^2 + \tilde{B}^1 \tilde{v}_{n+1}^1 \right] \end{aligned}\right\} \tag{78}$$

$$\therefore \delta = -\left\{ \sum_{i=1}^{S} L^i \begin{bmatrix} \overset{j=1}{\downarrow} & \cdots & \overset{j=dof(\lambda)}{\downarrow} \\ \dot{v}_{n+1}^i & \cdots & \dot{v}_{n+1}^i \end{bmatrix} \right\} \tag{79}$$

Notice that in order to compute $\delta$, we need the solution from step 1 at (n + 1) time instants. Hence for component sub-domains with $\xi^i > 1$, we only require solution at $\xi^i = \eta^i$. This information can be used to our advantage by saving only the results for $\eta^i = \xi^i$ in Eqs. (74). Matrix $\delta$ is similar to that obtained in the GC method for the case of $\Delta T = \Delta t$ [22]. As we can see, it may be calculated only once before entering the global time-stepping loop.

### 6.3. Step 3: Solve $Aw = f$

$$\begin{bmatrix} \tilde{M}^S & \cdots & 0 & 0 \\ \vdots & \ddots & \vdots & \vdots \\ 0 & \cdots & \tilde{M}^2 & 0 \\ 0 & \cdots & 0 & \tilde{M}^1 \end{bmatrix} \begin{bmatrix} \tilde{w}_{\xi^S}^S \\ \vdots \\ \tilde{w}_{\xi^2}^2 \\ \tilde{w}_{n+1}^1 \end{bmatrix} = \begin{bmatrix} \tilde{\tilde{F}}_{\xi^S}^S \\ \vdots \\ \tilde{\tilde{F}}_{\xi^2}^2 \\ \tilde{F}_{n+1}^1 - \tilde{N}^1 \tilde{U}_n^1 \end{bmatrix}$$

$$\tag{80}$$

Where: $\tilde{w}_{\xi^i}^i = \begin{bmatrix} \tilde{w}_1^i & \tilde{w}_2^i & \cdots & \tilde{w}_{\xi^i}^i \end{bmatrix}^T$

For $\Omega^1$ ($\xi^1 = 1$) we have:

$$\tilde{M}^1 \tilde{w}_{n+1}^1 = \tilde{F}_{n+1}^1 - \tilde{N}^1 \tilde{U}_n^1 \tag{81}$$

Since the above system of equations is solved under the action of external forces only, it is equivalent to solving the free problem for $\Omega^1$ as expressed by Eq. (36). For component sub-domains ($\xi^i > 1$), the solution is computed under the actions of external forces as well as known interface reactions ($\lambda_n$) and is updated after every intermediate time-step until $\xi^i = \eta^i$, that is until the solution is synchronized with (n+1) time instant. Step 3 represented as a sub-domain independent pseudo-code can be illustrated by Algorithm 3.

| for i = 1 to S | |
|---|---|
| for k = 1 to $\xi^i$ | |
| Input | |
| $M^i$ $\quad C^i \quad K^i$ | $F_k^i - (1 - k/\xi^i) {L^i}^T \lambda_{k-1}$ |
| $\ddot{w}_{k-1} \quad \dot{w}_{k-1} \quad w_{k-1}$ | $\Delta t^i \quad \beta^i \quad \gamma^i$ |
| Black-box 1: Newmark integration | |
| Output | |
| $\ddot{w}_k \quad \dot{w}_k \quad w_k \to \tilde{w}_k^i$ | |

Algorithm 3: Computing sub-domain response under the action of (applied) external forces and known interface reactions.

### 6.4. Step 4: Compute $y = \delta^{-1}(g - cw)$

$$cw = \begin{bmatrix} \tilde{B}^S & \cdots & \tilde{B}^2 & \tilde{B}^1 \end{bmatrix} \begin{bmatrix} \tilde{w}_{\xi^S}^S \\ \tilde{w}_{\xi^2}^2 \\ \vdots \\ \tilde{w}_{n+1}^1 \end{bmatrix} \tag{82}$$

$$\left. \begin{aligned} \therefore cw &= \tilde{B}^S \tilde{w}_{\xi^S}^S + \ldots + \tilde{B}^2 \tilde{w}_{\xi^2}^2 + \tilde{B}^1 \tilde{w}_{n+1}^1 \\ \therefore cw &= \tilde{B}^S \tilde{w}_{n+1}^S + \ldots + \tilde{B}^2 \tilde{w}_{n+1}^2 + \tilde{B}^1 \tilde{w}_{n+1}^1 \end{aligned} \right\} \tag{83}$$

$$\therefore cw = \sum_{i=1}^{S} L^i \dot{w}_{n+1}^i \tag{84}$$

$$\therefore y = \delta^{-1} \left\{ [0] - \sum_{i=1}^{S} L^i \dot{w}_{n+1}^i \right\} = -\delta^{-1} \sum_{i=1}^{S} L^i \dot{w}_{n+1}^i \tag{85}$$

$$\therefore y = \left\{ \sum_{i=1}^{S} L^i \begin{bmatrix} \overset{j=1}{\downarrow} & \cdots & \overset{j=dof(\lambda)}{\downarrow} \\ \dot{v}_{n+1}^i & \cdots & \dot{v}_{n+1}^i \end{bmatrix} \right\}^{-1} \sum_{i=1}^{S} L^i \dot{w}_{n+1}^i \tag{86}$$

Since we have already condensed out the interface reactions at intermediate time-steps, we only need to compute the Lagrange multipliers at time (T + ΔT). Equation (86) hence represents the direct solution of unknown Lagrange multipliers at time (T + ΔT).

### 6.5. Step 5: Compute x = w – vy

Global solution (combined from all component sub-domains) at T + ΔT and the local solution at intermediate time-steps ($\eta^i$ = 1, 2 …$\xi^i$) can now be computed as:

$$x = \begin{bmatrix} \tilde{\tilde{w}}^S_{\xi^S} \\ \vdots \\ \tilde{\tilde{w}}^2_{\xi^2} \\ \tilde{w}^1_{n+1} \end{bmatrix} - \begin{bmatrix} \tilde{\tilde{v}}^S_{\xi^S} \\ \vdots \\ \tilde{\tilde{v}}^2_{\xi^2} \\ \tilde{v}^1_{n+1} \end{bmatrix} y \tag{87}$$

Equation (87) is similar to Eq. (34) in the sense that the final solution is obtained as a sum of free and link responses, in this case however, one does not need to compute interface reactions at the intermediate time-steps. Now that the global solution is obtained for time (T + ΔT), we can repeat steps 3, 4 and 5 in order to compute the desired evolution of kinematic quantities. Aforementioned solution algorithm can be easily integrated into an existing structural finite element code. DD of the original problem may be performed prior to the analysis (user defined) or during analysis (based on a user defined criterion). Assuming user defined DD, we will now describe a general approach for MGMT implementation:

1. For every sub-domain $\Omega^i$ (i = 1,2 …S)
    a. Input nodes, coordinates, element connectivity's, material properties, boundary conditions, initial conditions, applied loads, time integration parameters and interface information.
    b. Form sub-domain arrays ($M^i$, $C^i$, $K^i$, $F^i$ etc.)
    c. Form Boolean projection matrices ($B^i$).
    d. Form interface constraint equations ($P^i$).
    e. Form interface connectivity matrices ($L^i = P^i B^i$).
    f. Form unit load response matrices using Algorithm 2, Step 1.
2. Compute interface condensation matrix (δ) using Step 2.
3. Factorize δ (this will be used in Step 4 to compute unknown Lagrange multipliers).
4. Compute initial accelerations based on initial conditions.
5. Initialize $\lambda_0 = 0.0$
6. Loop for total number of integration steps based on global time-step (ΔT).
    a. Computing sub-domain response under the action of (applied) external forces and known interface reactions using Algorithm 3, Step 3.
    b. Compute new interface reactions using Step 4.
    c. Compute global and local response using Step 5.
7. End loop.

### 7. Numerical examples

In this section we will present some examples that will help us evaluate the numerical stability, accuracy and efficiency of the MGMT method. We will analyze a 2D cantilever beam under forced vibrations; both transverse (Example 1) and longitudinal (Example 2). Domain under analysis and the transient forces applied at the free end are as shown in Fig. 8. In either case, the domain is discretized using 4-node quadrilateral elements (2 dof/node) with plane stress formulation, consistent mass matrix and zero damping. Isotropic linear elastic material properties with modulus of elasticity (E) = 2.07x10$^{11}$ N/m$^2$, Poisson's ratio (ν) = 0.3 and mass density (ρ) = 7.83x10$^3$ Kg/m$^3$ were used.

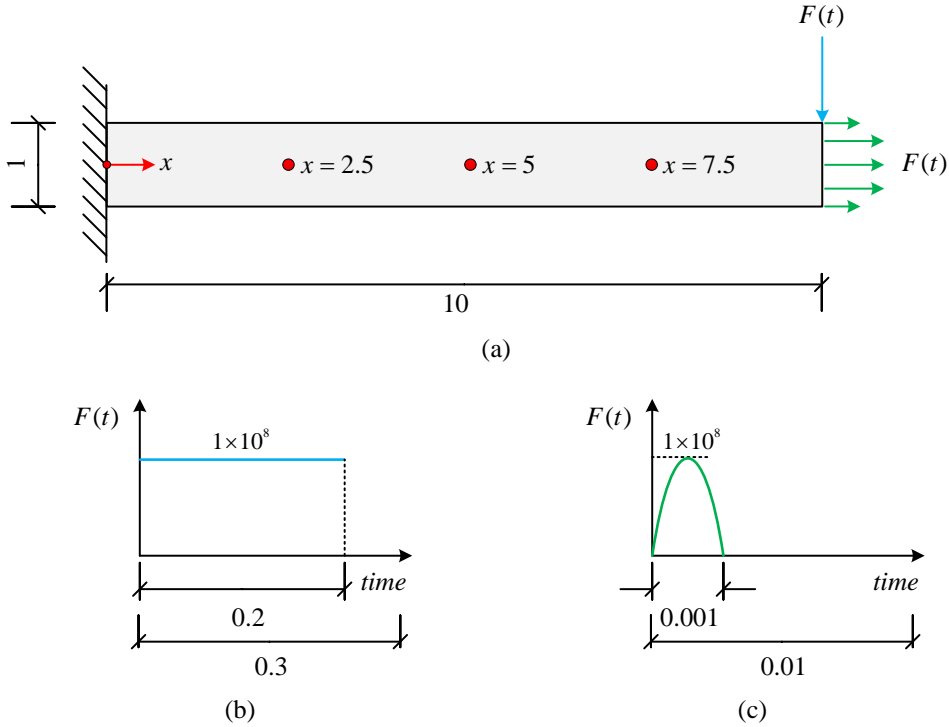

Fig. 8. (a) Domain under consideration (b) Transverse load (Example 1) (c) Longitudinal load (Example 2).

Different scenarios analyzed under these examples are as follows:

a. Case 1 – Reference cases. Uniform grid uniform time-step (UGUT) – Fig. 9 (a), Table 3.
b. Case 2 – Equally divided 4 sub-domains, I-I coupling (MGMT1) – Fig. 9 (b), Table 4.
c. Case 3 – Equally divided 2 sub-domains, I-E coupling (MGMT2) – Fig. 9 (c), Table 5.

Primary conclusions are drawn by comparing MGMT results with equivalent uniform grid uniform time-step (UGUT) cases. For MT coupling we will consider both, Implicit-Implicit and Implicit-Explicit algorithms. We will also compare the size of respective finite element problems (nodes, elements, number of equations, skyline storage requirements) and the required computation times (CPU solution time) to provide a basis for efficiency evaluation. When comparing results, we will use the root-mean-square error (RMSE) to provide a quantitative measure of accuracy between measured variables. This will enable us to quantify the error (on an average) between respective MGMT and UGUT cases. RMSE and Normalized RMSE (NRMSE) are computed as:

$$RMSE = \sqrt{\frac{1}{n}\sum_{i=1}^{n}\left(\text{var}_i^{MGMT} - \text{var}_i^{UGUT}\right)^2} \tag{88}$$

$$NRMSE = \frac{RMSE}{\text{var}_{max}^{UGUT} - \text{var}_{min}^{UGUT}} \times 100\% \tag{89}$$

Initial UGUT simulations were performed with increasing refinement, both spatial and temporal, in order to establish numerical convergence. Displacement at free end (x = 10) from Example 1 was measured and is plotted in Fig. 10. We can see that UGUT4 (Implicit) and UGUT5 (Explicit) yield converged results. Since these converged results provide best available approximation of the measured variable, primary comparisons are made against these cases. Since MGMT1 uses implicit algorithms in all its component sub-domains, it is compared v/s UGUT4 and since MGMT2 employs explicit algorithm in Ω2, it is compared v/s UGUT5. These pairs are highlighted in **bold** in forthcoming RMSE and NRMSE tables.

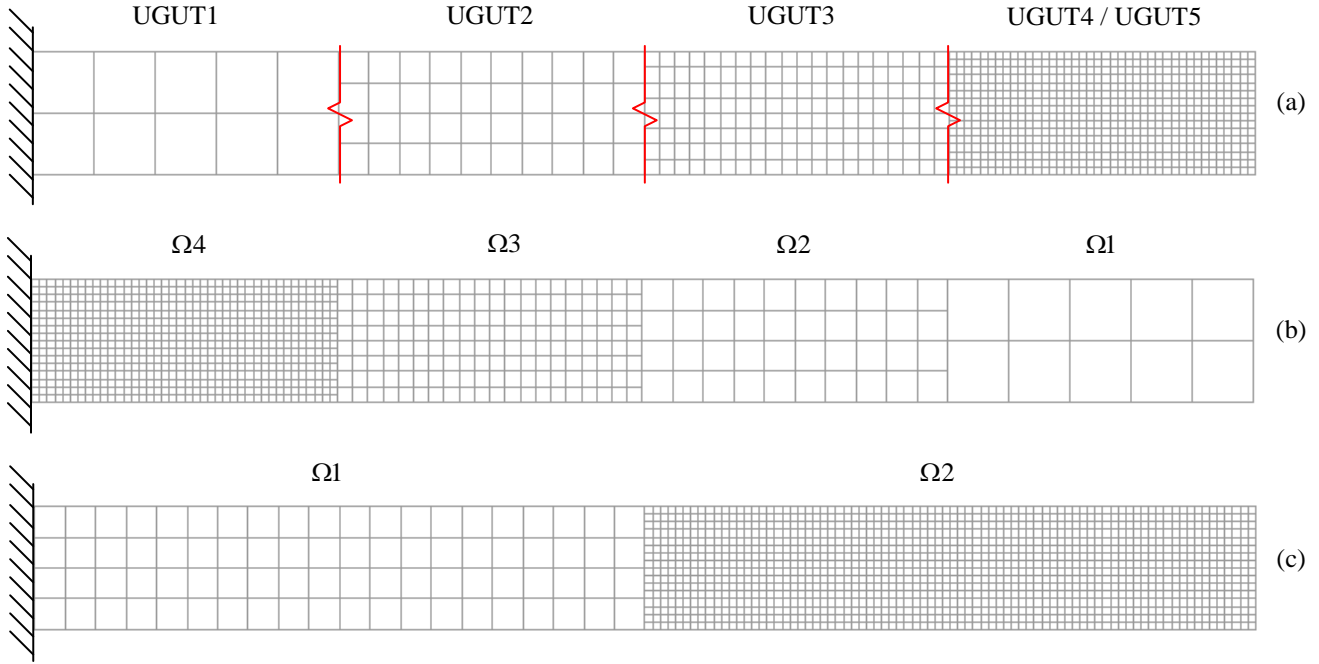

Fig. 9. Example grids – (a) Case 1 UGUT (b) Case 2 MGMT1 (c) Case 3 MGMT2.

| Sub-case ID | Grid spacing (H) | Newmark parameters | Time-step (ΔT) | |
|---|---|---|---|---|
| | | | Example 1 | Example 2 |
| UGUT1 | 0.5 | $\beta=0.25, \gamma=0.5$ (Implicit) | 1E-3 | 1E-4 |
| UGUT2 | 0.25 | $\beta=0.25, \gamma=0.5$ (Implicit) | 0.5E-3 | 0.5E-4 |
| UGUT3 | 0.125 | $\beta=0.25, \gamma=0.5$ (Implicit) | 0.25E-3 | 0.25E-4 |
| UGUT4 | 0.0625 | $\beta=0.25, \gamma=0.5$ (Implicit) | 0.125E-3 | 0.125E-4 |
| UGUT5 | 0.0625 | $\beta=0.0, \gamma=0.5$ (Explicit) | 1E-6 | 1E-6 |

Table 3: Case 1 – Reference cases (uniform grid uniform time-step).

| Sub-domain | Grid spacing (H) | Newmark parameters | Time-step (Δt) | | Time-step ratio ($\xi$) |
|---|---|---|---|---|---|
| | | | Example 1 | Example 2 | |
| $\Omega 1$ | 0.5 | $\beta=0.25, \gamma=0.5$ (Implicit) | 1E-3 | 1E-4 | 1 |
| $\Omega 2$ | 0.25 | $\beta=0.25, \gamma=0.5$ (Implicit) | 0.5E-3 | 0.5E-4 | 2 |
| $\Omega 3$ | 0.125 | $\beta=0.25, \gamma=0.5$ (Implicit) | 0.25E-3 | 0.25E-4 | 4 |
| $\Omega 4$ | 0.0625 | $\beta=0.25, \gamma=0.5$ (Implicit) | 0.125E-3 | 0.125E-4 | 8 |

Table 4: Case 2 – MGMT1 (4 sub-domains / dof ($\lambda$) = 65 / I-I coupling / Global time-step $\Delta T = 1E-3$).

| Sub-domain | Grid spacing (H) | Newmark Parameters | Time-step (Δt) | | Time-step ratio ($\xi$) | |
|---|---|---|---|---|---|---|
| | | | Example 1 | Example 2 | Example 1 | Example 2 |
| $\Omega 1$ | 0.25 | $\beta=0.25, \gamma=0.5$ (Implicit) | 0.5E-3 | 0.5E-4 | 1 | 1 |
| $\Omega 2$ | 0.0625 | $\beta=0.0, \gamma=0.5$ (Explicit) | 1E-6 | 1E-6 | 500 | 50 |

Table 5: Case 3 – MGMT2 (2 sub-domains / dof ($\lambda$) = 22 / I-E coupling / Global time-step $\Delta T = 0.5E-3$).

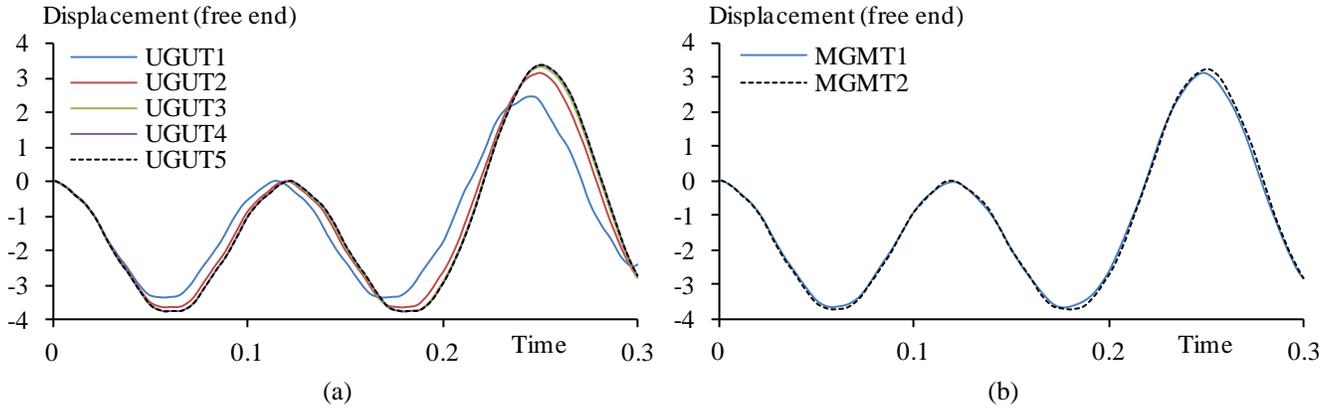

Fig. 10: Example 1 – displacement at free end (x = 10).

| v/s | UGUT1 | UGUT2 | UGUT3 | UGUT4 | UGUT5 |
|---|---|---|---|---|---|
| UGUT1 | - | 0.5667 (8.3%) | 0.7300 (10.34%) | 0.7758 (10.89%) | 0.7748 (10.87%) |
| UGUT2 | - | - | 0.0493 (0.7%) | 0.0509 (0.71%) | 0.0508 (0.71%) |
| UGUT3 | - | - | - | 0.0499 (0.7%) | 0.0482 (0.68%) |
| UGUT4 | Converged | | | | |
| UGUT5 | Converged | | | | |
| MGMT1 | 0.5272 (9.08%) | 0.0545 (0.8%) | 0.2112 (2.99%) | **0.2595 (3.64%)** | 0.2581 (3.62%) |
| MGMT2 | 0.6125 (10.56%) | 0.0509 (0.75%) | 0.1322 (1.87%) | 0.1811 (2.54%) | **0.1794 (2.52%)** |

Table 6: Example 1 – RMSE and NRMSE (%). Variable = displacement at free end.

Conformance between MGMT1 (v/s UGUT4) and MGMT2 (v/s UGUT5) is really good with less than 4% error. However in general, MGMT results are most in conformance with UGUT2, with less than 1% error. MGMT2, with loading applied in the finest discretization domain ($\Omega 2$), has relatively smaller errors compared to MGMT1, where loading is applied over the coarsest discretization ($\Omega 1$). Accordingly, when using MGMT method, it is recommended to use fine scale discretization in regions where externally applied forces are present or in sub-domains where local gradients or wave dynamics are to be captured accurately. The effect on results due to implicit v/s explicit algorithms is relatively small (0.02%) and can be seen in (MGMT2 v/s UGUT4) v/s (MGMT2 v/s UGUT5). Since the average error (over the entire range of reference results) is less than 4%, we can say that MGMT results are reasonably accurate and reliable.

In order to ensure numerical stability in MGMT coupling, we now look at global energies and augmented interface energy contributions from component sub-domains. We first compare global energies from UGUT4 v/s MGMT1, Fig. 11 with relative errors listed in Table 7, along with MGMT2 v/s UGUT5. Also listed is the RSME in interface energy indicating its mean variance about zero (since the corresponding UGUT interface energy is zero). Although the errors for kinetic and stiffness energy are relatively high, MGMT1 and MGMT2 satisfy global energy balance where kinetic + stiffness + interface energy = external energy. The residual after 0.2 sec (where the transient step load terminates) represents the steady state vibration where total energy oscillates between kinetic and stiffness (potential) energy. Augmented interface energy contributions from component sub-domains are plotted in Fig. 12 and Fig. 13 for MGMT1 and MGMT2 respectively. These plots clearly support Eq. (9) and Eq. (65) indicating the fact that sub-domain interface energies annihilate each other. Variance in interface energy, Table 7, is partly due to machine tolerance / round-off errors. Since it is orders of magnitude smaller than the total energy, its contribution can be considered to be negligible. By refining the increments used in numerical (trapezoidal) integration of Eqs. (12) and (13), errors can be significantly reduced. Another consequence of global stability requirement, continuity of velocities at sub-domain interfaces, can be verified and is plotted in Fig. 14 where we compare the vertical velocity at x = 5, as obtained from adjacent sub-domains of (a) MGMT1 and (b) MGMT2 in Example 1.

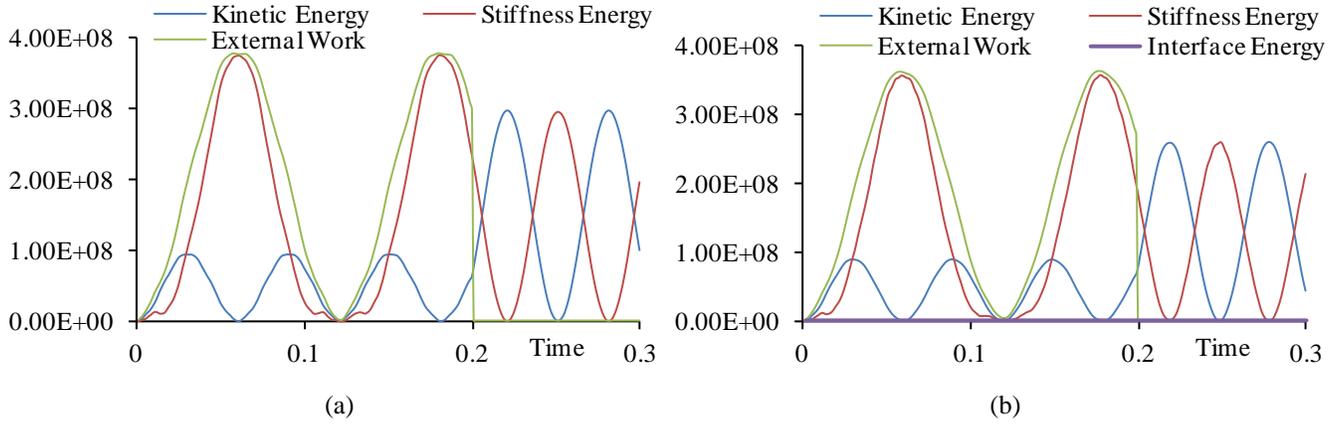

Fig. 11: Global energies (a) UGUT4 (b) MGMT1.

|  | Kinetic Energy | Stiffness Energy | Interface Energy |
|---|---|---|---|
| MGMT1 v/s UGUT4 | $2.37 \times 10^7$ (7.99%) | $2.54 \times 10^7$ (6.74%) | $6.26 \times 10^{-5}$ |
| MGMT2 v/s UGUT5 | $1.55 \times 10^7$ (5.22%) | $1.68 \times 10^7$ (4.47%) | $2.09 \times 10^{-5}$ |

Table 7: Example 1 – RMSE and NRMSE (%). Variable = Global energies.

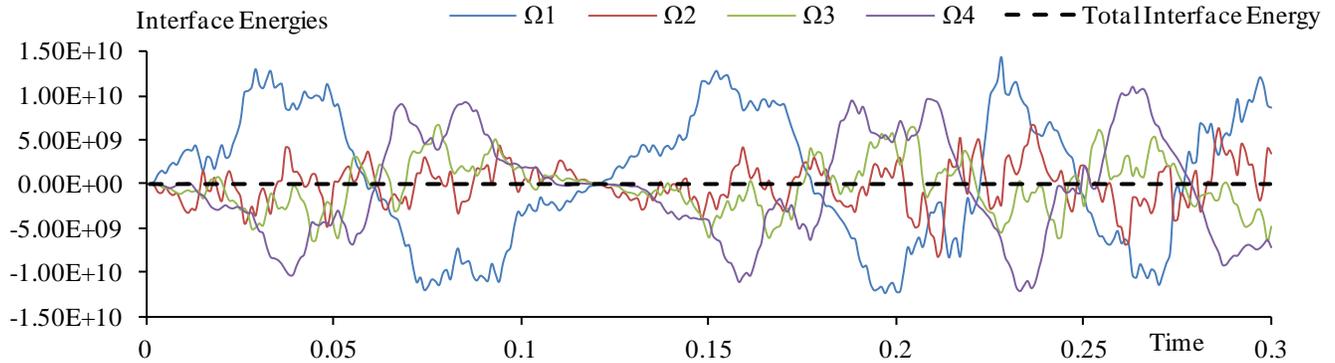

Fig. 12: Example 1 – MGMT1: Sub-domain interface energy contributions.

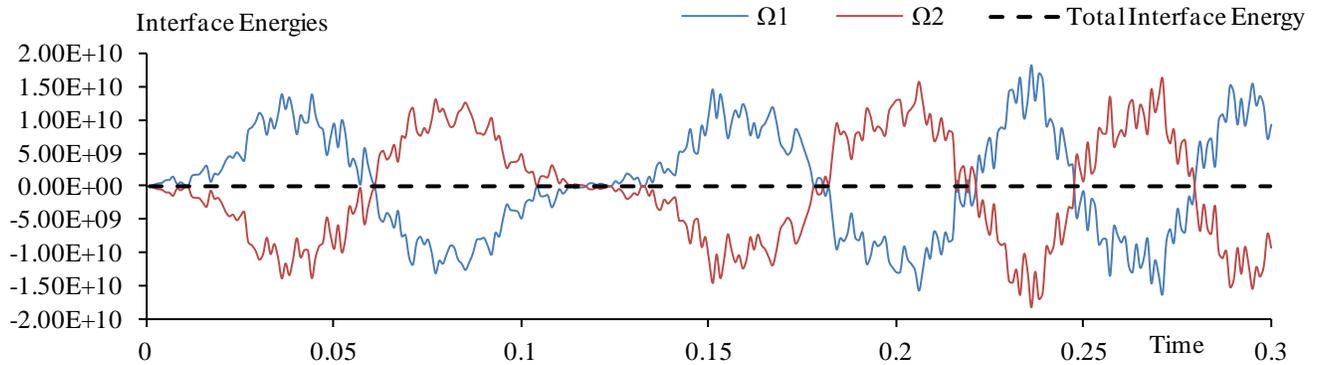

Fig. 13: Example 1 – MGMT2: Sub-domain interface energy contributions.

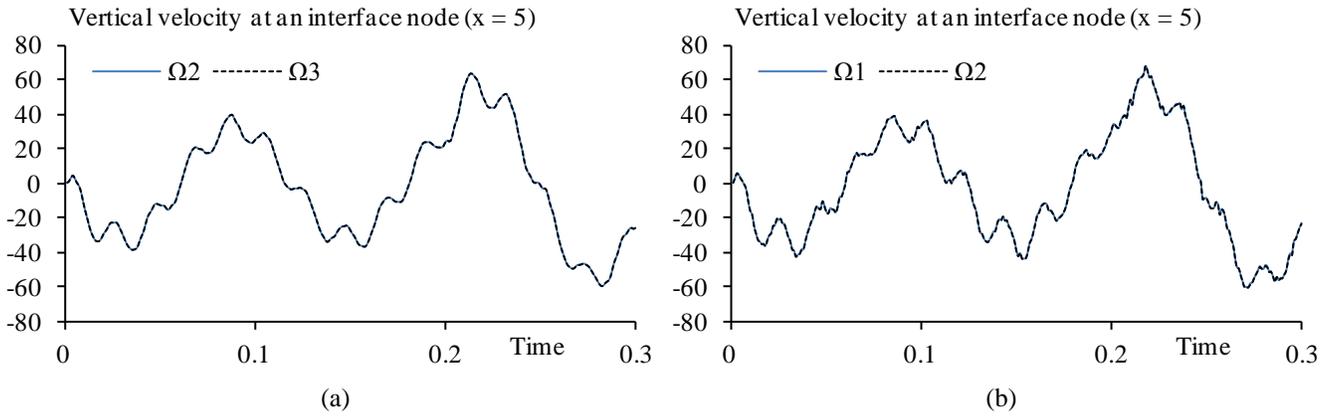

Fig. 14: Example 1 – Vertical velocity at an interface node as obtained from adjacent sub-domains (a) MGMT1 (b) MGMT2.

Primary focus in Example 2 is to evaluate the capability of MGMT method in simulating wave propagation across component sub-domains. Since structural waves are coupled in space and time, we need to ensure that MGMT formulation does not damage the characteristic features of a wave (such as its amplitude and phase) as it propagates through distinct spatial and temporal sub-domains. Fig. 15 shows the longitudinal stress (Sigma-xx) and Fig. 16 shows the longitudinal displacement wave as it propagates across the beam and measured at 4 different locations – x = 10, 7.5, 5 and 2.5. Stresses were computed at element level with 4 Gauss integration points and then projected onto sub-domain nodes. Table 9 and Table 10 list the RMSE and the corresponding normalized errors between UGUT and MGMT cases. For the longitudinal stress results, MGMT1 v/s UGUT4 at x = 10 ≡ Ω1 yields relatively high error since stresses in this sub-domain are computed over coarse grid. With a finer grid (MGMT2 v/s UGUT5 at x = 10 ≡ Ω2) the error sharply drops by almost 50%. Hence finely discretized sub-domains are recommended in regions where stress (or strain) gradients are critical so that the resulting errors can be kept to their minimum. The results for longitudinal displacement, Fig. 16, are very satisfactory with an average error of 1% (MGMT1 v/s UGUT4) and 2% (MGMT2 v/s UGUT5). Since the wave propagates seamlessly across MGMT sub-domains without any significant loss in amplitude or phase, we can consider these results to be reasonably accurate and reliable.

We now look at the principal advantage of using MGMT method as opposed to UGUT simulations by comparing the improvements in invested computational resources. Table 8 lists the total computation time (in sec) required to obtain the global solution of the problem, total number of nodes, elements, number of equations representing the primary dof and the total amount of skyline storage required for storing these variables. Total number of equations reduced by 65% (MGMT1 v/s UGUT4) and 45% (MGMT2 v/s UGUT5) and the corresponding simulation speedups are ~60% and ~40% respectively.

| Case ID | Nodes | Elements | Number of equations | Skyline storage | Solution time (sec) | |
|---|---|---|---|---|---|---|
| | | | | | Example 1 | Example 2 |
| UGUT1 | 63 | 40 | 120 | 1036 | 0.23400 | 0.07800 |
| UGUT2 | 205 | 160 | 400 | 5304 | 2.01241 | 0.65520 |
| UGUT3 | 729 | 640 | 1440 | 32132 | 22.23014 | 7.36324 |
| UGUT4 | 2737 | 2560 | 5440 | 221352 | 320.92525 | 107.53148 |
| UGUT5 | 2737 | 2560 | 5440 | 221352 | 24822.27071 | 830.14372 |
| | | | | | | |
| MGMT1 | 959 | 850 | 1884 | 74110 | 127.983220 | 41.808268 |
| MGMT2 | 1482 | 1360 | 2954 | 122615 | 13678.54938 | 490.279943 |

Table 8: Comparison of computational resources.

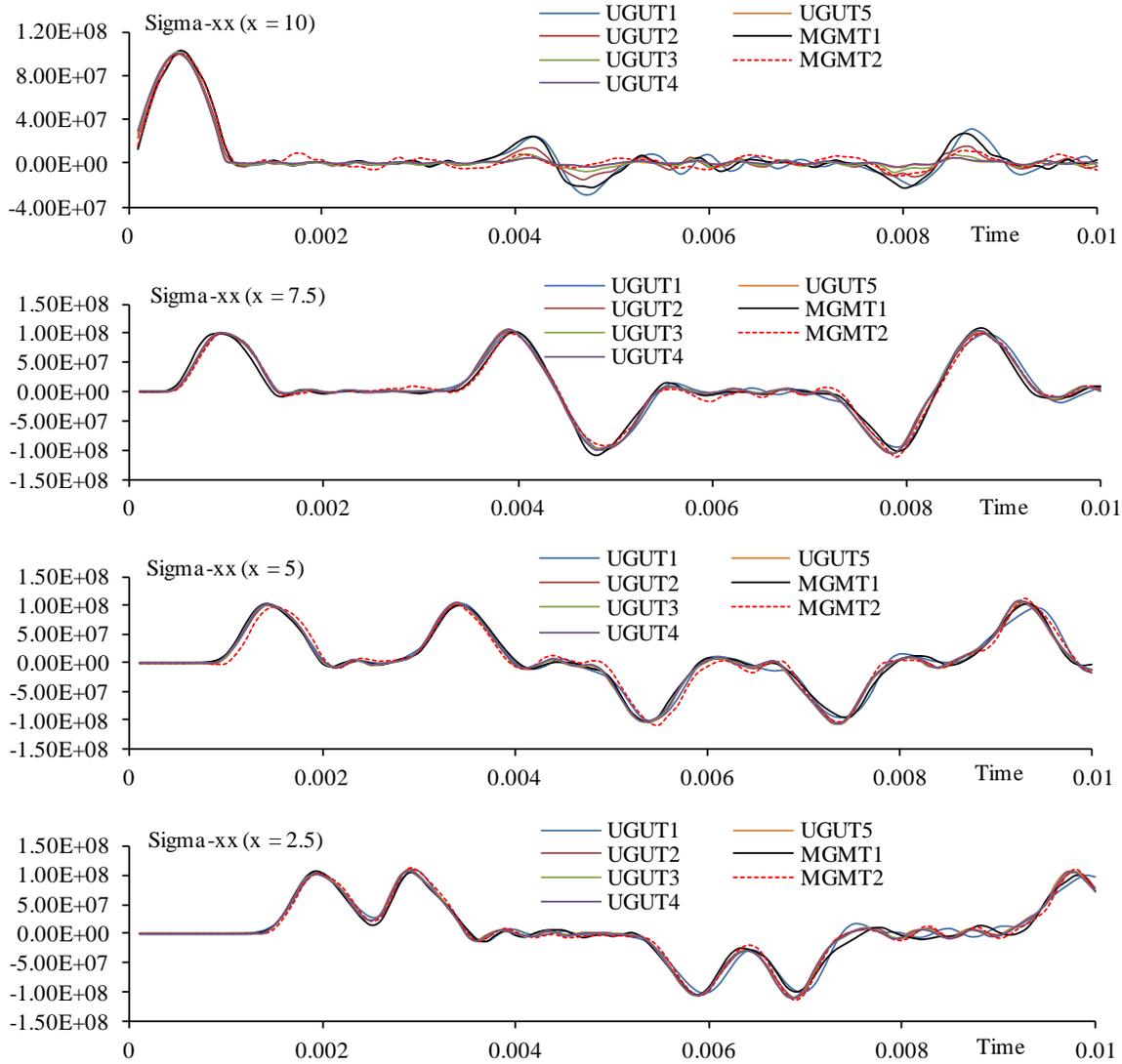

Fig. 15: Example 2 – Longitudinal stress (Sigma-xx) recorded at 4 different locations.

|        | UGUT1 | UGUT2 | UGUT3 | UGUT4 | UGUT5 |
|--------|-------|-------|-------|-------|-------|
| | | | x = 10 | | |
| MGMT1 | $3.78 \times 10^6$ (2.9%) | $4.89 \times 10^6$ (4.3%) | $7.19 \times 10^6$ (6.6%) | **$8.42 \times 10^6$ (8.0%)** | $8.48 \times 10^6$ (8.1%) |
| MGMT2 | $8.22 \times 10^6$ (6.4%) | $4.70 \times 10^6$ (4.1%) | $4.59 \times 10^6$ (4.2%) | $5.09 \times 10^6$ (4.9%) | **$5.06 \times 10^6$ (4.8%)** |
| | | | x = 7.5 | | |
| MGMT1 | $8.13 \times 10^6$ (4.0%) | $7.23 \times 10^6$ (3.5%) | $7.52 \times 10^6$ (3.6%) | **$7.60 \times 10^6$ (3.6%)** | $7.84 \times 10^6$ (3.7%) |
| MGMT2 | $7.40 \times 10^6$ (3.7%) | $6.38 \times 10^6$ (3.0%) | $6.84 \times 10^6$ (3.2%) | $6.93 \times 10^6$ (3.3%) | **$7.15 \times 10^6$ (3.4%)** |
| | | | x = 5.0 | | |
| MGMT1 | $7.08 \times 10^6$ (3.5%) | $5.06 \times 10^6$ (2.4%) | $5.73 \times 10^6$ (2.7%) | **$6.05 \times 10^6$ (2.9%)** | $6.11 \times 10^6$ (2.9%) |
| MGMT2 | $9.60 \times 10^6$ (4.7%) | $8.20 \times 10^6$ (3.8%) | $8.87 \times 10^6$ (4.2%) | $8.97 \times 10^6$ (4.3%) | **$9.27 \times 10^6$ (4.4%)** |
| | | | x = 2.5 | | |
| MGMT1 | $8.54 \times 10^6$ (4.1%) | $6.54 \times 10^6$ (2.9%) | $6.92 \times 10^6$ (3.2%) | **$7.10 \times 10^6$ (3.3%)** | $7.28 \times 10^6$ (3.4%) |
| MGMT2 | $8.33 \times 10^6$ (4.0%) | $4.41 \times 10^6$ (2.0%) | $5.29 \times 10^6$ (2.5%) | $5.51 \times 10^6$ (2.5%) | **$5.87 \times 10^6$ (2.7%)** |

Table 9: Example 2 – RMSE and NRMSE (%). Variable = Longitudinal stress (Sigma-xx).

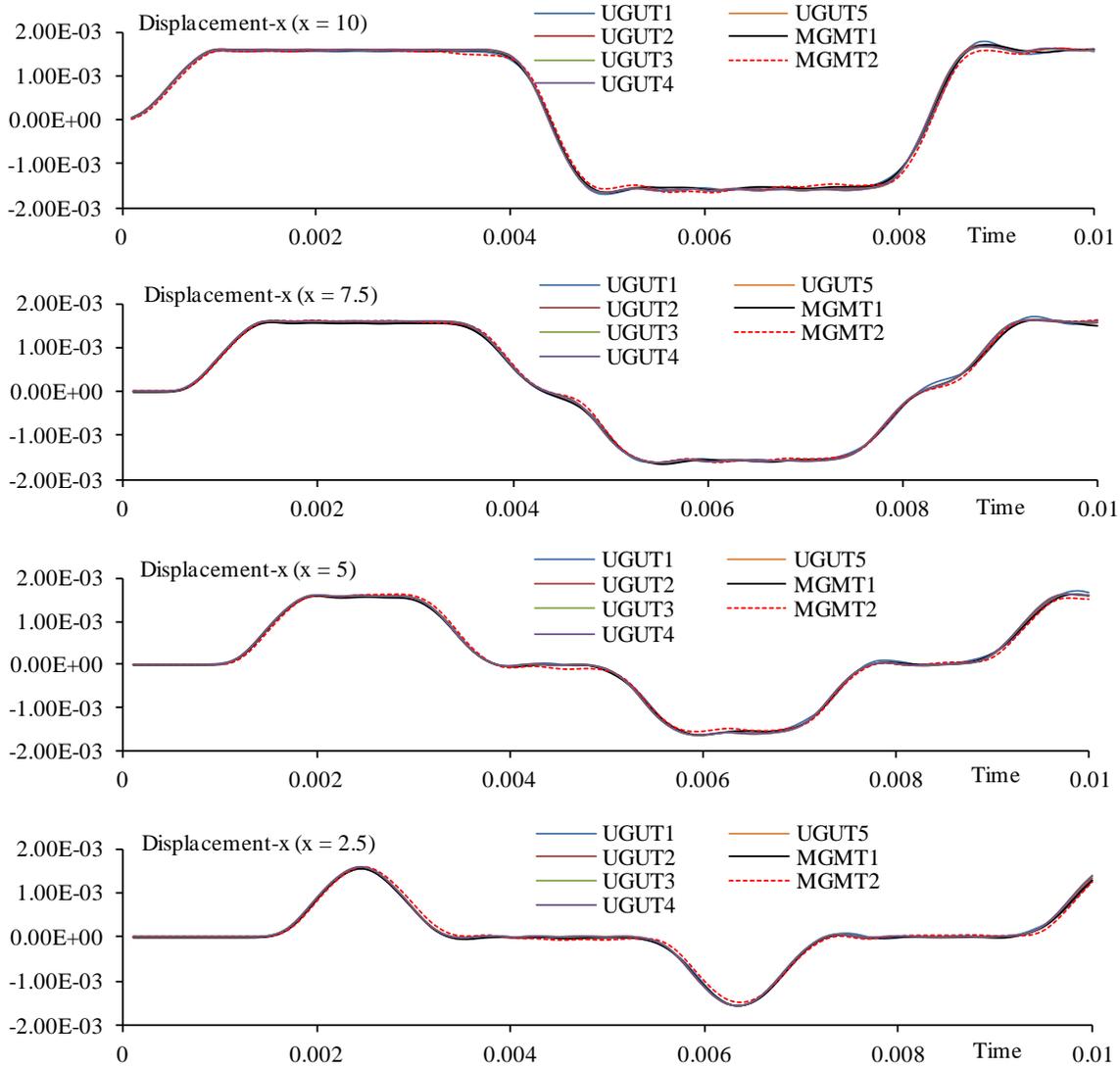

Fig. 16: Example 2 – Longitudinal displacement recorded at 4 different locations.

|  | UGUT1 | UGUT2 | UGUT3 | UGUT4 | UGUT5 |
|---|---|---|---|---|---|
|  | x = 10 | | | | |
| MGMT1 | $3.16 \times 10^{-5}$ (0.9%) | $2.82 \times 10^{-5}$ (0.8%) | $3.32 \times 10^{-5}$ (1.0%) | **$3.46 \times 10^{-5}$ (1.1%)** | $3.59 \times 10^{-5}$ (1.1%) |
| MGMT2 | $6.86 \times 10^{-5}$ (2.0%) | $6.82 \times 10^{-5}$ (2.0%) | $7.09 \times 10^{-5}$ (2.1%) | $7.15 \times 10^{-5}$ (2.2%) | **$7.25 \times 10^{-5}$ (2.2%)** |
|  | x = 7.5 | | | | |
| MGMT1 | $2.78 \times 10^{-5}$ (0.8%) | $2.34 \times 10^{-5}$ (0.7%) | $2.75 \times 10^{-5}$ (0.9%) | **$2.86 \times 10^{-5}$ (0.9%)** | $2.96 \times 10^{-5}$ (0.9%) |
| MGMT2 | $4.48 \times 10^{-5}$ (1.4%) | $3.76 \times 10^{-5}$ (1.2%) | $3.97 \times 10^{-5}$ (1.2%) | $4.03 \times 10^{-5}$ (1.2%) | **$4.11 \times 10^{-5}$ (1.2%)** |
|  | x = 5.0 | | | | |
| MGMT1 | $3.34 \times 10^{-5}$ (0.7%) | $3.26 \times 10^{-5}$ (0.6%) | $3.24 \times 10^{-5}$ (0.8%) | **$3.23 \times 10^{-5}$ (0.8%)** | $3.23 \times 10^{-5}$ (0.8%) |
| MGMT2 | $5.84 \times 10^{-5}$ (1.8%) | $5.78 \times 10^{-5}$ (1.8%) | $5.97 \times 10^{-5}$ (1.8%) | $6.01 \times 10^{-5}$ (1.8%) | **$6.08 \times 10^{-5}$ (1.9%)** |
|  | x = 2.5 | | | | |
| MGMT1 | $1.62 \times 10^{-5}$ (0.5%) | $1.46 \times 10^{-5}$ (0.5%) | $1.91 \times 10^{-5}$ (0.6%) | **$2.01 \times 10^{-5}$ (0.6%)** | $2.13 \times 10^{-5}$ (0.7%) |
| MGMT2 | $5.50 \times 10^{-5}$ (1.8%) | $6.17 \times 10^{-5}$ (2.0%) | $6.46 \times 10^{-5}$ (2.1%) | $6.52 \times 10^{-5}$ (2.1%) | **$6.60 \times 10^{-5}$ (2.1%)** |

Table 10: Example 2 – RMSE and NRMSE (%). Variable = Longitudinal displacement.

## 8. Conclusion

In this paper, we present a systematic approach to perform multiple grid and multiple time-scale (MGMT) simulations in linear structural dynamics. The formulation is largely based upon the fundamental principles of DDM, which allows us to selectively formulate, discretize and solve component sub-domains augmented with an interface condition. Coupled equations for decomposed sub-domains are linked together using Lagrange multipliers (used to represent interface reactions) and an appropriate interface condition is modeled by requiring that the interface energy (produced as a result of introducing interface reactions) is identically zero. We also show that enforcing this particular condition naturally results in the continuity of velocities across sub-domain interfaces. Conversely, we enforce continuity of velocities as an interface condition and show, using Energy method, that the resulting interface energy contributions are zero. Subsequently, we show that as long as the stability requirements are satisfied within the time integration of component sub-domains, MGMT coupling is stable and energy preserving. We also present a step-by-step approach for connecting multiple grids (conforming or non-conforming) using the fundamental principles of Mortar FEM used to derive interface constraint equations. Multiple time-scale coupling is achieved using Newmark family of algorithms and an energy preserving approach is used for the condensation of intermediate interface reactions. Final set of equations for MGMT sub-domains are solved using block elimination and Crout factorization and we present the corresponding algorithm to obtain the global solution at synchronous time-steps. Its implementation is discussed in details and in a fashion that allows easy integration into existing FE codes.

Through a series of numerical examples and comparisons made against uniform grid uniform time-scale (UGUT) simulations, we show that MGMT method is stable, efficient and reasonably accurate in modeling concurrent multi-domain simulations. Normalized RMSE are within acceptable range, ensuring a good conformance between MGMT and converged UGUT cases. Stability was verified numerically by ensuring zero interface energy (RMSE is in the order of $10^{-5}$) and continuity of velocities across sub-domain interfaces. We also show that the characteristic features of a structural wave (stress and displacement) are preserved as it propagates seamlessly across MGMT sub-domains without any significant change in its amplitude or phase.

Finally, we look at the computational resources by comparing the total number of equations and the corresponding skyline storage requirements for every case. We clearly see that MGMT cases require very few degrees of freedom in comparison to UGUT4/5 (converged uniform grid uniform time-scale, I/E). The loss in available degrees of freedom is clearly reflected in the NRMSE; however it is relatively insignificant compared to the advantage gained in simulation speedup.

**Acknowledgements**

We would like to thank Dr. Kunik Lee and FHWA for their continual support. This research is sponsored under DOT Cooperative Agreement No. DTFH61-10-00005